\documentclass[12pt]{article}
\usepackage{amssymb}
\def\theequation{\thesection.\arabic{equation}}

\newtheorem{theorem}{Theorem}

\def\thetheorem{\thesection.\arabic{theorem}}

\def\theprop{\thesection.\arabic{prop}}
\newtheorem{lemma}[theorem]{Lemma}
\def\thelemma{\thesection.\arabic{lemma}}
\newtheorem{cor}[theorem]{Corollary}
\def\thecor{\thesection.\arabic{cor}}

\def\theexam{\thesection.\arabic{exam}}

\def\theremark{\thesection.\arabic{remark}}
\newcommand{\eqa}{\begin{eqnarray}}
\newcommand{\eeqa}{\end{eqnarray}}
\newcommand{\beq}{\begin{equation}}
\newcommand{\eeq}{\end{equation}}
\newcommand{\nn}{\nonumber}

\newcommand{\pal}{\partial}

\newcommand{\al}{\alpha}

\newcommand{\ve}{\epsilon}

\newcommand{\pf}{\noindent{\it Proof \ }}
\newcommand{\tr}{{\rm tr}}
\newcommand{\res}{{\rm res}}
\newcommand{\lgl}{\log{L}}
\newcommand{\ld}{\Lambda}
\newcommand{\epf}{$\quad$\hfill
\raisebox{0.11truecm}{\fbox{}}\par\vskip0.4truecm}

\begin{document}

\title {\LARGE Virasoro Symmetries of the Extended Toda  Hierarchy}
\author{ {Boris Dubrovin${}^*$ \ \ Youjin Zhang${}^{**}$}\\
{\small ${}^{*}$ \ SISSA, Via Beirut 2--4, 34014 Trieste, Italy}\\
{\small ${}^{**}$ \ Department of Mathematical Sciences, Tsinghua
University, Beijing 100084, P.R.China }} \maketitle
\begin{abstract}
We prove that the extended Toda hierarchy of \cite{CDZ} admits nonabelian Lie
algebra of infinitesimal symmetries isomorphic to the 
half of the Virasoro algebra. The 
generators $L_m$, $m\geq -1$ of the Lie algebra 
act by linear differential operators onto the tau
function of the hierarchy. 
We also prove that the tau function of a generic solution 
to the extended Toda hierarchy
is annihilated by a combination of the Virasoro operators and the flows of the
hierarchy.
As an application we show 
that the validity of the
Virasoro constraints for the $CP^1$ Gromov-Witten invariants and their
descendents
implies that their generating function is the logarithm of a particular 
tau function of the
extended Toda hierarchy.
\end{abstract}

\setcounter{equation}{0} \setcounter{theorem}{0}

\section{Introduction}
The extended Toda hierarchy was introduced in \cite{Z, getzler, CDZ} in an 
attempt to encode the
recursion relations among the $CP^1$ Gromov-Witten invariants into that of a 
hierarchy of
integrable systems. As it was shown in \cite{CDZ}, this hierarchy can be 
represented
in a Lax pair formalism through the Lax operator
\beq
L=\Lambda+v(x)+e^{u(x)} \Lambda^{-1}.
\eeq
The functions $v, u$  serve as the dependent variables for the hierarchy 
with spatial
variable $x$ and $\Lambda=e^{\ve\pal_x}$ is the shift operator, $\ve$ is a small
parameter. We will also introduce a two-component vector 
$$
w=(w^1, w^2), \quad w^1 := v, \, w^2:= u
$$
to use it in the formulae where many summations enter.

The flows of the hierarchy
are defined via Lax representation
\beq\label{td-flow-Lax}
\ve\frac{\pal L}{\pal t^{\beta,q}}=[A_{\beta,q},L] :=A_{\beta,q} L-L
A_{\beta,q},\quad \beta=1,2;\ q\ge 0.
\eeq
Here the operators $A_{\beta,q}$  have the expression
\beq\label{def-a-1}
A_{1,q}=\frac2{q!}\left[L^q (\lgl-c_q)\right]_+,\quad
A_{2,q}=\frac1{(q+1)!}\left[L^{q+1}\right]_+,
\eeq
\beq\label{def-c}
c_q=1+\frac12+\dots+\frac1{q}
\eeq
with the positive part $B_+$ of a difference operator $B=\sum B_k \ld^k$ given 
by $B_+=\sum_{k\ge 0} B_k
\ld^k$. The logarithm of the operator $L$ is defined as follows \cite{CDZ}. Let
us first introduce the
dressing operators $P$ and $Q$ of the form
\beq\label{PQ-2}
P=\sum_{k\ge 0} p_k \ld^{-k},\quad Q=\sum_{k\ge 0} q_k
\ld^{k},\quad p_0=1
\eeq
such that
\beq\label{PQ-1}
L=P\Lambda P^{-1}=Q\Lambda^{-1} Q^{-1}.
\eeq
Then we define
\beq\label{def-log}
\lgl:=\frac12 \left(P\ve\pal_x P^{-1}-Q\ve \pal_x Q^{-1}\right).
\eeq

The equations (\ref{td-flow-Lax}) for $\beta=2$ coincide with the standard
flows of the Toda hierarchy. In particular for $q=0$ one obtains the equations
of Toda lattice
\eqa\label{toda-stand}
&&
{\pal v\over \pal t^{2,0}} = {1\over \ve} \left( e^{u(x+\ve) }-
e^{u(x)}\right)=\sum_{k\geq 0} {\ve^{k}\over (k+1)!} \pal_x^{k+1} e^{u}
\nn\\
&&
{\pal u\over \pal t^{2,0}} = {1\over \ve} \left( v(x) - v(x-\ve)\right)
=\sum_{k\geq 0} (-1)^{k} {\ve^k\over (k+1)!} \pal_x^{k+1} v
\eeqa
written in the interpolated form. To return to the original discrete setup 
of \cite{Toda} one
is to introduce the lattice variables
$$
u_n:= u(n\ve), ~~v_n:= v(n\ve).
$$
The parameter $\ve$ plays the role of the mesh
of the lattice.
Another part, for $\beta=1$ is a new one. For $q=0$ one obtains just the spatial
translations
\beq\label{spat}
{\pal v\over \pal t^{1,0}} = {\pal v\over \pal x}, \quad
{\pal u\over \pal t^{1,0}} = {\pal u\over \pal x}.
\eeq
The flow for $\beta=1$ and $q=1$ is less trivial:
\eqa\label{t-1-1}
&& 
\frac{\pal v}{\pal t^{1,1}}=v v_x+\frac1{\ve}\left[e^{u(x+\ve)}
\left({\cal B}_- u(x+\ve)-2\right)
-e^{u(x)}\left({\cal B}_- u(x-\ve)-2\right)\right],\nn\\
&&
\frac{\pal u}{\pal t^{1,1}}=\frac1{\ve} 
\left[v(x) \left({\cal B}_- u(x)-2\right)-v(x-\ve)
\left({\cal B}_- u(x-\ve)-2\right)\right.\nn\\
&&
\quad \quad \quad \quad\quad \left.
+{\cal B}_+ v(x+\ve)-{\cal B}_+ v(x-\ve)\right]
\eeqa
where the operators ${\cal B}_\pm$ are defined by
\eqa
&&{\cal B}_+  :=(\Lambda-1)^{-1}\ve \pal_x  =\sum_{k\ge 0} 
\frac{B_k}{k!}(\ve \pal_x)^k  ,\nn\\
&&
{\cal B}_-  :=(1-\Lambda^{-1})^{-1}\ve \pal_x  =\sum_{k\ge 0} 
\frac{B_k}{k!}(-\ve \pal_x)^k  .\label{def-calB}
\eeqa
Here $B_k$ are the Bernoulli numbers.

The flows of the extended Toda hierarchy can be represented 
as Hamiltonian systems
\beq\label{HF}
\frac{\pal w^\al}{\pal t^{\beta,q}}=\{ w^\al (x), H_{\beta,q}\}_1\equiv
U_1^{\al\gamma}\,
\frac{\delta H_{\beta,q}}{\delta w^\gamma (x)}
\eeq
with the Hamiltonian operators
\beq\label{fpb}
U_1^{11}=U_1^{22}=0,\quad U_1^{12}=\frac1{\ve} (\Lambda-1),\quad U_1^{21}=\frac1{\ve} (1-\Lambda^{-1}),
\eeq
and the Hamiltonians
$$
H_{\beta,q}=\int h_{\beta,q} dx,\quad \beta=1,2; \ q\ge -1
$$ 
with the densities $h_{\beta,q}$ defined by
\eqa\label{def-h}
&&
h_{1,q}=\frac2{(q+1)!}\,\res\left[ L^{q+1}
(\lgl-c_{q+1})\right]=\res \,A_{1,q+1},
\nn\\
&&
h_{2,q}=\frac1{(q+2)!}\res\,  L^{q+2} =\res \,A_{2,q+1}.
\eeqa
By definition the residue of a difference operator
$$
A=\sum_{k\in {\mathbb Z}} a_k \Lambda^k
$$
is given by
$$
\res\, A := a_0.
$$
Note that 
\beq\label{casi}
h_{1,-1} = {\cal B}_-u(x), \quad h_{2,-1} = v(x)
\eeq
are densities of Casimirs of the Poisson bracket, i.e.
$$
\{~.~, H_{1,-1}\}_1 = \{~.~, H_{2,-1}\}_1 \equiv 0.
$$

Denote ${\cal
A}$ the graded ring of formal power series of the form 
$ f=\sum_{k\ge 0} f_k(w,w_x,\dots) \ve^k $,
where $f_k$ are polynomials of $v, u, e^{\pm u},v^{(m)}, u^{(m)},\ m\ge 1$. 
The gradation 
is defined by
\beq
\label{grad-A} \deg v^{(m)}=1-m,\ \deg u^{(m)}=-m,\ {\rm for} \ m\geq 0, \ \deg
e^{u}=2, \ \deg \ve=1.
\eeq
As it was shown in \cite{CDZ}, all the Hamiltonian densities (\ref{def-h})
as well as the right hand sides of the flows of the extended Toda hierarchy
are homogeneous elements of the ring ${\cal A}$.

Introduce functions $\Omega_{\al,p;\beta,q}\in {\cal A}$ by
\beq\label{omega-def-1} \frac1{\ve}\,
(\ld-1)\Omega_{\al,p;\beta,q}:=\frac{\pal h_{\al,p-1}}{\pal
t^{\beta,q}} =\left\{\begin{array}{ll} \frac2{p!}\,
\res\left(\left[A_{\beta,q},L^{p}
(\lgl-c_{p})\right]\right),& \al=1,\\
\frac1{(p+1)!}\,\res \left[A_{\beta,q},L^{p+1}\right],& \al=2.
\end{array}\right.
\eeq
Existence of such functions and an important property of $\tau$-{\it symmetry}
$$
\Omega_{\beta,q;\al, p}=\Omega_{\al,p;\beta,q}
$$
was established in \cite{CDZ}.
These elements of the ring ${\cal A}$  
are uniquely determined by the above formulae and 
by the homogeneity condition
$$
\deg \Omega_{\al,p;\beta,q}=p+q+1+\mu_\al+\mu_\beta.
$$
Here
$$
\mu_1= -{1\over 2}, \quad \mu_2 = {1\over 2}.
$$

In this paper we will consider the solutions to the extended Toda
hierarchy in the class of formal series in $\ve$
\eqa\label{sol-eps}
&&
w^\alpha(x, {\bf t};\ve) = \sum_{k\geq 0} \ve^k w^\al_k(x, {\bf t}), 
\quad \alpha=1,
\, 2,
\\
&&
{\bf t}= (t^{1,0}, t^{2,0}, t^{1,1}, t^{2,1}, \dots)
\nn
\eeqa
As it follows from the definition (\ref{def-a-1}),
$$
A_{1,0}=\pal_x
$$
modulo terms commuting with $L$. So 
$$
{\pal\over \pal t^{1,0}} = {\pal\over\pal x},
$$
i.e., the solution depends on $x$, $t^{1,0}$ only via the combination $x+t^{1,0}$.
We will therefore often suppress the explicit dependence on $x$ in the formulae.

\vskip 0.4truecm
\noindent{\bf Definition}(\cite{CDZ})\ For any solution $v(x,{\bf t};\ve)$,
$u(x, {\bf t};\ve)$
of
the extended Toda hierarchy there exists a function 
$$
\tau=\tau(x, {\bf t};\ve)=e^{\sum_{g\geq 0} \ve^{2g-2} {\cal F}_g(x,{\bf t})}
$$ 
such
that the functions $\Omega_{\al,p; \beta,q}$ evaluated on this solution can be
represented in the form
\beq\label{def-omega}
\Omega_{\al,p;\beta,q}=\ve^2\frac{\pal^2\log\tau} {\pal t^{\al,p}\pal t^{\beta,q}}
\eeq
for any $\al,\beta=1,2,\ p,q\ge 0$.
It is called the {\it tau function} of the solution (\ref{sol-eps})
of the extended Toda hierarchy.

In particular, we have
\beq\label{td-20}
h_{\al,p}=\ve (\ld-1)\frac{\pal\log\tau}{\pal t^{\al,p+1}},\quad \al=1,2,\ 
p\ge -1
\eeq
\eqa
&&
v=\ve (\Lambda-1) \frac{\pal\log\tau}{\pal t^{2,0}} = \ve {\pal\over\pal
t^{2,0}}
\log {\tau(x+\ve, {\bf t}; \ve)\over \tau(x, {\bf t}; \ve)},
\nn\\
&&
u=(\Lambda-1)(1-\Lambda^{-1}) \log\tau = \log {\tau(x+\ve, {\bf t}; \ve) 
\tau(x-\ve, {\bf t}; \ve)\over \tau^2(x, {\bf t}; \ve)}.\label{td-21}
\eeqa

Another important property of this hierarchy is that, apart from its Hamiltonian
structure described above, it also possesses a second Hamiltonian structure which is
compatible with the first one (see (\ref{secpb}) below). 
The bihamilonian structure and the tau symmetry property of the
extended Toda hierarchy imply, due to a general theorem of \cite{DZ3},
{\it quasi-triviality}
of the extended Toda hierarchy. The precise formulation of the quasitriviality
property in the case of interest is 
given by the following theorem:

\begin{theorem}[\cite{DZ3}]\label{thint}
Any solution $v, u$ of the extended Toda hierarchy is obtained from a solution $v_0, u_0$ of
the dispersionless extended Toda hierarchy through the quasi-Miura transformation of the form
\eqa
&&v=v_0+\sum_{g\ge 1} \ve^{2g-1}(\Lambda-1)
\frac{\pal F_g(w_0,\dots,w_0^{(3g-2)})}{\pal t^{2,0}},\nn\\
&&u=u_0+\sum_{g\ge 1} \ve^{2g-2}(\Lambda-1)(1-\Lambda^{-1})F_g(w_0,\dots,w_0^{(3g-2)}),
\eeqa
and the corresponding tau function of the solution admits the following 
{\rm genus expansion}
\beq\label{genus-exp}
\log\tau=\ve^{-2}\,\log\tau^{[0]}+\sum_{g\ge 1} \ve^{2g-2} F_g(w_0,\dots,w_0^{(3g-2)}).
\eeq
Here $\tau^{[0]}=\tau^{[0]}(x, {\bf t})$ is the tau function for the solution 
$v_0, u_0$ of the {\it dispersionless} extended Toda hierarchy, i.e., it is related
to the leading term of the solution (\ref{sol-eps}) 
$$
w_0= (v_0, u_0)
$$
via
\eqa\label{genus-exp0}
&&
v_0(x,{\bf t}) = {\pal^2 \log \tau^{[0]}(x, {\bf t})\over \pal x \pal t^{2,0}}
\nn\\
&&
u_0(x,{\bf t}) = {\pal^2 \log \tau^{[0]}(x, {\bf t})\over \pal x^2 }.
\eeqa
\end{theorem}

Recall that the dispersionless extended Toda hierarchy is obtained from 
the extended Toda hierarchy by setting $\ve=0$. All the flows
of the dispersionless
extended Toda hierarchy are systems of hydrodynamic type, i.e. systems
of two first order quasilinear PDEs. For $\beta=1$, $q=0$ the dispersionless
flow still coincides with the spatial translations (\ref{spat}). For $\beta=2$,
$q=0$ one obtains
\eqa
&&
{\pal v\over \pal t^{2,0}} = {\pal\over \pal x} e^u
\nn\\
&&
{\pal u\over \pal t^{2,0}} = {\pal v\over \pal x}.
\nn
\eeqa
Eliminating $v$ yields the so-called long wave limit of the Toda lattice
equations
$$
u_{tt} = \left( e^u \right)_{xx}
$$
where $t=t^{2,0}$. The dispersionless limit of (\ref{t-1-1}) reads
\eqa
&&
{\pal v\over \pal t^{1,1}} ={\pal\over \pal x} \left[ \frac12 v^2 +(u-1)
e^u\right]
\nn\\
&&
{\pal u \over \pal t^{1,1}} = {\pal \over \pal x} (u\, v).
\nn
\eeqa
Changing the sign of time $t=-t^{1,1}$ one identifies these with the equations
of motion of one-dimensional polytropic gas with the speed $v$ and density $u$
and the equation of state of the form
$
p=(u^2-2u+2)e^u -2
$.

It is time to remind to the reader that the theory of dispersionless (extended)
Toda hierarchy can be nicely encoded \cite{D0, D1}
in terms of a particular two-dimensional
Frobenius manifold $M_{\rm Toda}$. The latter can be identified with the quantum cohomology
of complex projective line 
$$
M_{\rm Toda}=QH^* (CP^1).
$$
Alternatively, the Frobenius manifold in question is isomorphic to the orbit
space of the simplest extended affine Weyl group \cite{compo}
$$
M_{\rm Toda}= {\mathbb C}^2 /\tilde W(A_1).
$$
Denote $v$, $u$ the coordinates on $M_{\rm Toda}$. The potential of the Frobenius manifold
reads
\beq\label{frob}
F= \frac12 v^2 u +e^u.
\eeq
The third derivatives of the potential define multiplication law of tangent
vectors at each point of $M$ such that $\pal / \pal v$ is the unity and
$$
{\pal \over \pal u} \cdot {\pal \over \pal u}= e^u {\pal \over \pal v}.
$$
The flat metric on $M$ reads
\beq\label{metric}
<~,~>=2 du\, dv.
\eeq
The Euler vector field is
$$
E=v {\pal \over  \pal v}  +2 {\pal \over \pal u}.
$$
We will not remind here the universal construction, due to \cite{npb}
of a ``dispersionless''
integrable hierarchy valid for an arbitrary Frobenius  manifold $M$. The hierarchy
can be considered as an infinite family of pairwise commuting flows on the
loop space ${\cal L}(M)$. All these flows are represented by first order
quasilinear PDEs; for this reason this hierarchy is called dispersionless.
The word ``hierarchy'' means that the systems of integrable PDEs are organized
by means of the action of a bihamiltonian recursion operator.

In \cite{DZ3} we addressed the problem of extending the correspondence
$$
{\rm Frobenius}~ {\rm manifolds} ~\to~ {\rm hierarchies} ~ {\rm of}~ 
{\rm integrable  } ~ {\rm PDEs}
$$
to an arbitrary Frobenius manifold. We proved that, indeed such a {\it universal
correspondence} exists for an arbitrary semisimple $M$ provided a suitable
completion of the loop space ${\cal L}(M)$ is made 
allowing to work with infinite order PDEs. By definition
semisimplicity of a Frobenius manifold means that, for a generic point $w\in M$
the algebra on $T_wM$ is semisimple. We leave as an exercise to the reader
to verify that $M_{\rm Toda}$ is a semisimple Frobenius manifold. So, according to
the main result of \cite{DZ3}, there exists an integrable hierarchy associated
with the Frobenius manifold $M_{\rm Toda}$. The main aim of the present paper is
to identify this integrable hierarchy with the extended Toda hierarchy
(\ref{td-flow-Lax}).

The crucial role in such identification, apart from the Lax representation and
tau structure obtained in \cite{CDZ}, play {\it Virasoro symmetries}. According
to our paper \cite{DZ2} for an arbitrary Frobenius manifold there exists a
universal construction of second order linear differential operators 
$$
L_m=L_m(\ve^{-1}{\bf t},
\ve{\pal\over\pal {\bf t}}), \, m\ge -1
$$ 
satisfying the Virasoro commutation relations
\beq\label{comm-L}
[L_i, L_j]=(i-j) L_{i+j}, \quad i, \, j \geq -1.
\eeq
For $M_{\rm Toda}$ these Virasoro operators are
given by the following explicit expressions (cf. \cite{eguchi})
\eqa
&& L_{-1}=\sum_{k\geq 1} {t}^{\alpha,k} {\pal\over
\pal t^{\alpha, k-1}}+{1\over \epsilon^2 } {t}^{1,0} {t}^{2,0},
\nn\\
&& L_0 =\sum_{k\geq 1} k\, \left( {t}^{1,k} {\pal\over \pal
t^{1,k}}+{t}^{2,k-1} {\pal\over \pal t^{2,k-1}}\right)+2\sum_{k\ge 1}
{t}^{1,k} \frac{\pal} {\pal t^{2,k-1}}+\frac1{\ve^2} \left(
{t}^{1,0}\right)^2,
\nn\\
&& L_m={\epsilon^2} \sum_{k=1}^{m-1} k!\, (m-k)!
{\pal^2\over \pal t^{2,k-1} \pal t^{2,k-m-1}} +\sum_{k\geq 1}
{(m+k)!\over (k-1)!}\left(
{t}^{1,k} {\pal\over \pal t^{1, m+k}}\right.\nn\\
&&\left.\quad + {t}^{2,k-1}{\pal\over \pal t^{2,k-1}}\right)
+2\,\sum_{k\geq 0} \alpha_m(k) {t}^{1,k} {\pal \over \pal {t}^{2,
m+k-1}}, \quad m\ge 1.\label{def-vir-td}
\eeqa
where
$$
\alpha_m(0)= m!, ~~\alpha_m(k) = {(m+k)!\over (k-1)!} \sum_{j=k}^{m+k} {1\over
j}, ~k>0.
$$
Again, we will not reproduce here the {\it universal construction} of the
Virasoro operators; however, we will give below the free field realization of
these operator for the $M_{\rm Toda}$ case obtained in our paper \cite{DZ3}. 

It is the {\it Virasoro invariance} property of a hierarchy of integrable PDEs that allows to reconstruct
it uniquely starting from a given semisimple Frobenius manifold, 
along with bihamiltonian
structure and existence of a tau function. That means that, the
{\it linear} action of the Virasoro operators onto tau functions 
define symmetries of the hierarchy\footnote{In \cite{DZ3} we called this
property {\it linearization of the Virasoro symmetries}. The reason for this
name was the following one. The dispersionless hierarchy on ${\cal L}(M)$
constructed in \cite{npb} is always invariant with respect to an action of half
of the Virasoro algebra \cite{DZ2}. However, the generators of the Virasoro action do not
act linearly onto dispersionless tau function $\tau^{[0]}({\bf t})$. The full
hierarchy constructed in \cite{DZ3} is a deformation of the dispersionless one
that transforms the nonlinear action of the Virasoro symmetries to the linear
one.}.  

The central result of this paper is the following\newline
\vskip -0.3truecm
\noindent {\bf Main Theorem}\, {\em 1. The transformations
\beq\label{vir-sym}
\tau\mapsto \tau + \delta \, L_m \tau , \quad \delta \to 0
\eeq
for any $m\geq -1$ are {\rm infinitesimal symmetries} of the extended Toda
hierarchy, i.e., given a tau function $\tau=\tau({\bf t}; \ve)$ of a solution
$v({\bf t};\ve)$, $u({\bf t}; \ve)$
of the form (\ref{sol-eps}), the functions 
\eqa
&&
\tilde v({\bf t};\ve)=v({\bf t};\ve)+\ve \delta (\Lambda-1) {\pal\over \pal t^{2,0}} {L_m\tau\over
\tau} +O(\delta^2)
\nn\\
&&
\tilde u({\bf t};\ve)=u({\bf t};\ve) +\delta (\Lambda-1)(1-\Lambda^{-1})
{L_m\tau\over\tau}+O(\delta^2)
\nn
\eeqa
satisfy equations of the extended Toda hierarchy modulo terms of order
$O(\delta^2)$.

\noindent 2. For a generic solution (\ref{sol-eps}) of the extended Toda hierarchy, 
the corresponding tau function
satisfies the Virasoro constraints
\beq
{L}_m(\ve^{-1}({\bf t}-{\bf c}(\ve)), \ve{\pal\over \pal {\bf t}}) \tau=0,\quad m\ge -1
\eeq
for some ${\bf c}(\ve)$.
Here ${\bf c}(\ve)= \left( c^{\al,p}(\ve)\right)$ is a collection of formal power
series in $\ve$.
}
\vskip -0.2 truecm

We will describe the class of generic solutions in Section 4. Note that in the
above formulae we omit writing explicitly the $x$-dependence since $x$ enters
only through the combination $x+t^{1,0}$.

The above two theorems and the uniqueness result of \cite{DZ3} imply 

\begin{cor} \label{cor-01} The hierarchy of PDEs associated, 
according to \cite{DZ3}, with
the Frobenius manifold $M_{\rm Toda}$ coincides with (\ref{td-flow-Lax}).
\end{cor}

According to \cite{DZ3}, the functions $F_g=F_g(w_0; \dots,
w_0^{(3g-2)})$ in the genus expansion 
(\ref{genus-exp}) are uniquely determined
by the loop equation that was introduced in \cite{DZ3} for any semisimple 
Frobenius 
manifold. For $M_{\rm Toda}$ the loop equation will be discussed in the Section
5. It will also be explained how one can compute Gromov - Witten invariants
of $CP^1$ and their descendents of any genus using the loop equation.

Due to the uniqueness of solution
of the loop equation and the validity of the Virasoro constraints for $CP^1$ 
Gromov-Witten invariants \cite{givental2}, we
have the following

\begin{cor}\label{cor-00} The generating function of the $CP^1$ Gromov-Witten 
invariants and their descendents is uniquely determined by the system of
Virasoro constraints. It coincides with the logarithm of the tau function
of a particular solution of the extended Toda hierarchy.
\end{cor}

This particular solution for the $CP^1$ Gromov-Witten invariants will  
be described in Section 5.

Summarizing, we can say that, in the theory of Gromov - Witten invariants
of $CP^1$
and their descendents 
the extended Toda hierarchy (\ref{td-flow-Lax}) plays the role similar to one
played by the KdV hierarchy in the Kontsevich - Witten formulation
\cite{Witten1, Kon, Witten2}
of the intersection theory on the Deligne - Mumford moduli spaces of stable punctured
curves.
\vskip 0.5truecm
\noindent{\bf Acknowledgments.} The researches of B.D. were
partially supported by Italian Ministry of Education research grant Cofin2001
``Geometry
of Integrable Systems''. The researches of Y.Z. were partially supported
by the Chinese National Science Fund for Distinguished Young Scholars grant
No.10025101 and the Special Funds of Chinese Major Basic Research Project
``Nonlinear Sciences''. Y.Z. thanks 
Abdus Salam International Centre for Theoretical Physics and SISSA
where part of the work was done for the hospitality.

\setcounter{equation}{0} \setcounter{theorem}{0}

\section{Some formulae for the functions $\Omega_{\al,p;\beta,q}$}

We present in this section some important identities for the functions
$\Omega_{\al,p;\beta,q}$ defined in(\ref{omega-def-1}) 
and for the Hamiltonians $H_{\al,p}$.

Let us first recall the definitions of \cite{CDZ} of the logarithmic
$x$-derivatives of the dressing operators $P$ and $Q$. Let us look for these
logarithmic derivatives in the form
\eqa\label{log-der}
&&
\ve P_x P^{-1} =\sum_{j\geq 1} b_j \Lambda^{-j}
\nn\\
&&
\ve Q_x Q^{-1} =\sum_{j\geq 0} c_j \Lambda^{j}.
\eeqa

\begin{lemma} \label{old-lemma} There exist unique elements $b_j, \, c_j\in{\cal A}$ homogeneous
of the degrees
$$
\deg b_j = j, \quad \deg c_j=-j
$$
such that the
operators (\ref{log-der}) satisfy the following system of equations
$$
\res\,\left( [\ve P_x P^{-1}, L^m]-\ve \pal_x L^m \right) =
\res\,\left( [\ve Q_x Q^{-1}, L^m]-\ve \pal_x L^m \right) =0, \quad m\geq 1.
$$
\end{lemma}

Proof see in \cite{CDZ}. For example,
$$
b_1=-{\cal B}_+v, \quad c_0={\cal B}_- u, \quad c_1=e^{-u(x+\ve)}{\cal B}_-v.
$$

\begin{lemma}
The following identities hold true
\eqa
&&\frac{\pal\Omega_{\al,p;\beta,q}}{\pal
v}=\Omega_{\al,p-1;\beta,q}+
\Omega_{\al,p;\beta,q-1}+(\delta_{\al,1}\delta_{\beta,2}+\delta_{\al,2}
\delta_{\beta,1}) \delta_{p,0}\delta_{q,0};\label{omega-du1}\\
&&\left(\sum_{m\ge 0}v^{(m)}\frac{\pal}{\pal v^{(m)}}+2
\frac{\pal}{\pal u}
\right) \Omega_{\al,p;\beta,q}\nn\\
&&=(p+q+1+\mu_\al+\mu_\beta)
\Omega_{\al,p;\beta,q}+R^{\gamma}_\al\Omega_{\gamma,p-1;\beta,q}+
R^{\gamma}_\beta\Omega_{\al,p;\gamma,q-1}\nn\\
&&\quad +
2\,\delta_{\al,1}\delta_{\beta,1}\delta_{p,0}\delta_{q,0}.
\label{omega-dE1}
\eeqa
Here the numbers $\mu_\al$ and a $2\times 2$ matrix $R=\left(
R^\gamma_\beta\right)$ are defined by
\beq\label{def-mu-R}
\mu_1=-\mu_2=-\frac12,\quad R^\gamma_\beta=2 \delta^\gamma_2 \delta_{\beta,1}.
\eeq
\end{lemma}
\pf Let us consider the difference operators
\eqa
&&
B:= {\pal\over\pal v} P\ve \pal_x P^{-1} =- \sum_{j\geq 1} {\pal b_j\over \pal
v}\Lambda^{-j}
\nn\\
&&
C:= {\pal\over\pal v} Q\ve \pal_x Q^{-1} =- \sum_{j\geq 0} {\pal c_j\over \pal
v}\Lambda^{j}
\nn
\eeqa
The coefficients $b_j\in{\cal A}$, $c_j\in{\cal A}$ were defined in
(\ref{log-der}).
We want to show that 
\beq\label{tozh}
LB=1, \quad LC=-1.
\eeq
Indeed, differentiating the identity
$$
[P\ve \pal_x P^{-1}, L]=0
$$
with respect to $v$
we derive that the operators $B$ and $L$ commute. Therefore
$$
[B-L^{-1},L^m]=0, \quad m\geq 1.
$$
Since $b_1=-{\cal B}_+v$, 
$$
{\pal b_1\over \pal v}=-1.
$$
So the expansion of the difference operator $B-L^{-1}$ begins with
$\Lambda^{-2}$. From the equations
$$ 
\res [B-L^{-1},L^m]=0 \quad {\rm for}~{\rm any} ~m\geq 1
$$
we prove that the coefficient of $\Lambda^{-2}$ of the operator $B-L^{-1}$ is a
constant.
Since
$$
\deg {\pal b_j\over \pal v} =j-1,
$$
the degree of this coefficient, as an element of ${\cal A}$, is equal to 1. So,
the coefficient must be equal to zero. Continuing this process we prove that
$$
B=L^{-1} =\Lambda^{-1} +O(\Lambda^{-2}).
$$ 
In a similar way we can prove that 
$$
C = - L^{-1} =-e^{-u(x+\ve)} \Lambda +O(\Lambda^2).
$$

Now, using the definition (\ref{def-log}) of $\log L$ we obtain
$$
{\pal\over\pal v} \log L = \frac12 (B-C) = L^{-1}.
$$
Therefore
\beq\label{toda-H3}
L^q \frac{\pal \lgl}{\pal v}=L^{q-1}.
\eeq
From the last equation it readily follows that
\eqa
&&\frac{\pal}{\pal v}\left[\frac2{p!} L^p
(\lgl-c_p)\right]
=\frac2{(p-1)!} L^{p-1} (\lgl-c_{p-1}),\nn\\
&& \frac{\pal}{\pal v}\left[\frac1{(p+1)!} L^{p+1}\right]
=\frac1{p!} L^{p}.\nn
\eeqa
These two equations yield, together with
$$
h_{1,-1}=(1-\ld^{-1})^{-1}\ve\pal_x u,\quad h_{2,-1}=v
$$
and the definition (\ref{omega-def-1}),
the formula (\ref{omega-du1}).

To prove the second formula of the Lemma, let us introduce
the operators ${\cal H}$ and ${\cal E}$ that act on the space of
difference operators of the form $\sum b_k \ld^k$ with
coefficients in ${\cal A}$ as follows:
\eqa
&&{\cal H} \sum b_k \ld^k=\sum k b_k \ld^k,\\
&&{\cal E}\sum b_k \ld^k=\sum \left(\sum_{m\ge 0} v^{(m)}
\frac{\pal b_k}{\pal v^{(m)}}+2 \frac{\pal b_k}{\pal u}\right)
\ld^k.
\eeqa
It is easy to see that these operators satisfy
the Leibnitz rule
\beq
{\cal P}\left[(\sum a_k \ld^k)(\sum b_l
\ld^l)\right] =\left({\cal P}\sum a_k \ld^k\right)\sum b_l \ld^l+
\sum a_k \ld^k \left({\cal P}\sum b_l \ld^l)\right).
\eeq
Here ${\cal P}={\cal H}$ or ${\cal P}={\cal E}$. By our definition it
follows that
\beq\label{H-E-l}
\left({\cal H}+{\cal E}\right) L=L.
\eeq
Due to Lemma \ref{old-lemma}
\beq\label{H-E-lgl}
\left({\cal H}+{\cal E}\right)\lgl=1.
\eeq
Now by using (\ref{H-E-l}) and (\ref{H-E-lgl}) we obtain \eqa
&&\frac1{\ve}(\ld-1) {\cal E}\Omega_{2,p;1,q}
=\frac1{(p+1)!}\,{\cal E}\res \left[A_{1,q},L^{p+1}\right]\nn\\
&&=\frac1{(p+1)!}\,\left({\cal H}+{\cal E}\right)
\res \left[A_{1,q},L^{p+1}\right]\nn\\
&&=\frac1{(p+1)!}\,\res \left[q A_{1,q}+2
A_{2,q-1},L^{p+1}\right]+
\frac1{p!}\,\res \left[A_{1,q},L^{p+1}\right]\nn\\
&&=(p+q+1)\frac1{\ve}(\ld-1)
\Omega_{2,p;1,q}+2\,\frac1{\ve}(\ld-1) \Omega_{2,p;2,q-1}. \eeqa
So, from the homogeneity condition for $\Omega_{\al,p;\beta,q}$ we arrive at
$$
{\cal E}\Omega_{2,p;1,q}=(p+q+1)\Omega_{2,p;1,q} + 2 \Omega_{2,p;2,q-1}.
$$
Other cases of the formula (\ref{omega-dE1}) can be proved in a
similar way. The Lemma is proved. \epf

\begin{lemma}\label{prop-1}
The variational derivatives of the Hamiltonians $H_{\beta,q}$ are given by the
following formulae
\eqa &&\frac{\delta
H_{\beta,q}}{\delta v (x)}=h_{\beta,q-1}=\ve (\ld-1) \frac{
\pal\log\tau}{\pal t^{\beta,q}},\label{delta-H-u1}\\
&&\frac{\delta H_{\beta,q}}{\delta u(x)}=\Omega_{2,0;\beta,q}=\ve^2
\frac{ \pal^2\log\tau}{\pal t^{2,0}\pal
t^{\beta,q}}.\label{delta-H-u2} \eeqa
\end{lemma}

\pf From the Hamiltonian representation (\ref{HF}) of the extended Toda
hierarchy we have
\beq
\frac{\pal v}{\pal t^{\beta,q}}=\frac1{\ve} \left(\Lambda-1\right)\frac{\delta
H_{\beta,q}}{\delta u},\quad
\frac{\pal u}{\pal t^{\beta,q}}=\frac1{\ve}
\left(1-\Lambda^{-1}\right)\frac{\delta H_{\beta,q}}{\delta v}.
\eeq
On the other hand, the formulae (\ref{td-20}) and (\ref{td-21}) imply that
\beq
\frac{\pal v}{\pal t^{\beta,q}}=\frac1{\ve}\left(\Lambda-1\right)
\Omega_{2,0;\beta,q},\quad
\frac{\pal u}{\pal t^{\beta,q}}=\frac1{\ve}\left(1-\Lambda^{-1}\right)
h_{\beta,q}.
\eeq
So, due to the homogeneity property of the densities of the Hamiltonians we
arrive at the expressions
(\ref{delta-H-u1}) and (\ref{delta-H-u2}).
The Lemma is proved.\epf

\setcounter{equation}{0} \setcounter{theorem}{0}

\section{Virasoro operators for the Frobenius manifold $M_{\rm Toda}$}

Recall that, according to the general algorithm of \cite{DZ3} the Virasoro
operators associated with a given Frobenius manifold are obtained by the following
free field realization\footnote{In \cite{DZ2} we used a different free field
realization of the Virasoro operators inspired by \cite{eguchi2}.}
 that we now present for the example of $M_{\rm Toda}$. Let
$a_{1,p}$, $a_{2,p}$, $p\in {\mathbb Z}$ be free bosonic operators satisfying
the following commutation relation\footnote{Note change of notations: 
in \cite{DZ3} we used half integer
labels.}
\beq\label{komm}
[a_{1,p}, a_{2,q}]= (-1)^p \delta_{p+q+1,0}.
\eeq
Introduce vectors of operators
\beq
{\bf a}_p:= (a_{1,p}, a_{2,p}), \quad p\in {\mathbb Z}.
\eeq
Consider the generating function
\beq\label{g-f}
{\bf f}(z):= \sum_{p\in {\mathbb Z}} {\bf a}_p z^{p+\mu} z^R
=\sum_{p\in {\mathbb Z}} {\bf a}_p \left( \matrix{ 
z^{p-\frac12} & 0 \cr
2 z^{p+\frac12} \log z & z^{p+\frac12} \cr}\right).
\eeq
Here the diagonal matrix $\mu$ and nilpotent matrix $R$ correspond to the {\it
spectrum at the origin} of the Frobenius manifold $M_{\rm Toda}$
\beq\label{mu-r}
\mu=\left( \matrix{ -\frac12 & 0 \cr & \cr 0 & \frac12\cr}\right), \quad
R=\left( \matrix{0 & 0\cr 2 & 0\cr}\right).
\eeq
The current $\phi^{(\nu)}(\lambda)$ for any non-integer $\nu$ is defined by
a (suitably defined) Laplace-type transform
\eqa\label{lap}
&&
\phi^{(\nu)}(\lambda)= \int_0^\infty {dz\over z^{\frac12 +\nu}} \, e^{-\lambda\, z} {\bf f}(z)
\nn\\
&&
= \sum_{p\in {\mathbb Z}} {{\bf a}_p\over \lambda^{p-\nu}} 
\left( \matrix{ \Gamma(p-\nu) & 0 \cr {2\over \lambda} \left[
\Gamma'(p-\nu+1)-\log\lambda\, \Gamma(p-\nu+1)\right] & \frac1{\lambda} 
\Gamma(p-\nu+1)\cr}\right).
\eeqa
The generating function of the regularized Virasoro operators $L_m^{(\nu)}$
is defined by the following quadratic combination of the derivatives of the
currents
\eqa\label{stress}
&&
T^{(\nu)}(\lambda) =\sum_{m\in{\mathbb Z}} {L^{(\nu)}_m\over \lambda^{m+2}}
\nn\\
&&
=
\frac12 :\pal_\lambda \phi^{(-\nu)} G(\nu) \left[ \pal_\lambda \phi^{(\nu)}\right]^T:
+{1\over 4 \lambda^2} \tr \left( {1\over 4} - \mu^2\right)
\nn\\
&&
=\sum_{p, q\in{\mathbb Z}}{1\over \lambda^{p+q+3}} :{\bf a}_p 
M_{pq}(\nu,\lambda)
\, {\bf a}_q^T :,\quad {\rm where}
\\
&&
M_{pq}(\nu,\lambda) = \left(\matrix{ 0 &  {\sin \pi \nu\over 2\pi}
\Gamma(p+\nu+1)\Gamma(q-\nu+2)\cr
 & \cr
-{\sin \pi \nu\over 2\pi}
\Gamma(q-\nu+1)\Gamma(p+\nu+2) & \frac1{\lambda}\pal_\nu \left[
{\sin\pi \nu\over \pi} \Gamma(p+\nu+2)\Gamma(q-\nu+2)\right]\cr}\right).
\nn
\eeqa
In this formula the Gram matrix $G(\nu)$ reads
\beq\label{gram-g}
G(\nu) ={1\over 2\pi} \left[ e^{\pi\, i\, R}e^{\pi\, i (\mu+\nu)} + 
e^{-\pi\, i\, R}e^{-\pi\, i (\mu+\nu)}\right] \eta^{-1}
={1\over \pi} \left( \matrix{ 0 & \sin\pi\nu\cr
-\sin\pi\nu & 2\pi \cos \pi\nu\cr}\right),
\eeq
$\eta$ is the Gram matrix of the metric (\ref{metric}), i.e., in our case
$$
\eta=\eta^{-1}=\left( \matrix{ 0 & 1\cr 1 & 0\cr}\right);
$$
the term with the trace in (\ref{stress}) does not contribute since $\mu^2 =
{1\over 4}$. The normal ordering in (\ref{stress}) is defined in such a way
that all the operators $a_{1,p}$ and $a_{2,p}$ with nonnegative $p$ are to be
written on the right. So
\eqa
&&
L_m^{(\nu)} =\sum_{p+q=m-2} :a_{2,p} a_{2,q}: \pal_\nu \left[ {\sin\pi\nu\over \pi}
\Gamma(p+\nu+2) \Gamma(q-\nu+2)\right].
\nn\\
&&
+\sum_{p+q=m-1} :a_{1,p} a_{2,q}: {\sin\pi\nu\over 2\pi}\left[
\Gamma(p+\nu+1)\Gamma(q-\nu+2)-\Gamma(p-\nu+1)\Gamma(q+\nu+2)\right].
\nn
\eeqa
Using the standard identities of the theory of gamma-functions
$$
\Gamma(x+1) = x\, \Gamma(x), \quad \Gamma(x)\Gamma(1-x) ={\pi\over \sin\pi x}
$$
we finally obtain the following expression for the regularized Virasoro
operators
\eqa\label{vir-oper}
&&
L_m^{(\nu)}=\frac12 \sum_p (-1)^{p+1} \left[
{\Gamma(\nu-p+m+1)\over \Gamma(\nu-p)}+ {\Gamma(-\nu-p+m+1)\over \Gamma(-\nu-p)}
\right]
:a_{1,p}a_{2, m-p-1}: 
\nn\\
&&
+\sum_p (-1)^p \pal_\nu \left[ {\Gamma(m-\nu-p)\over
\Gamma(-\nu-p-1)}\right] \, : a_{2,p} a_{2,m-p-2}:.
\eeqa
The operators
$$
L_{-1}^{(\nu)}= \sum_p (-1)^{p+1} :a_{1,p} a_{2,-p-2}:
$$
and
$$
L_0^{(\nu)} = \sum_p (-1)^p p\, :a_{1,p} a_{2, -p-1}: +\sum_p (-1)^{p+1}
:a_{2,p} a_{2,-p-2}:
$$
do not depend on $\nu$. The operators (\ref{vir-oper}) with $m>0$ depend
polynomially on $\nu$; those with $m<-1$ depend rationally on $\nu$. Therefore
there exist limits
$$
L_m:= \lim_{\nu\to 0} L_m^{(\nu)}, \quad m\geq -1.
$$

To arrive at the Virasoro operators given above in the Introduction one is 
to use the following realization of the bosonic operators $a_{\al,p}$
\beq\label{bose}
a_{\al,p} = \left\{ \matrix{ \ve {\pal\over \pal t^{\al,p}}, & p\geq 0\cr
 & \cr
\ve^{-1} (-1)^{p+1} t_{\al,-p-1}, & p<0\cr}.\right.
\eeq
Here we use the matrix $\eta$ for lowering the indices
$$
t_{\alpha,p}:= \eta_{\al\beta}t^{\beta,p}.
$$

In the next Section we will prove that the linear action of the 
Virasoro operators $L_m$ with $m\geq -1$ defines infinitesimal symmetries
of the extended Toda hierarchy.

\setcounter{equation}{0} \setcounter{theorem}{0}

\section{Proof of the main theorem}
We first consider the following system of Euler-Lagrange equations:
\eqa\label{EL-toda}
&&
\sum_{p\ge 0} {\tilde t}^{\al,p}\,\frac{\delta
H_{\al,p-1}}{\delta v(x)}=0
\nn\\
&&
\sum_{p\ge 0} {\tilde t}^{\al,p}\,\frac{\delta
H_{\al,p-1}}{\delta u(x)}=0
\eeqa
where
\beq {\tilde t}^{\al,p}=t^{\al,p}-c^{\al,p}(\ve)+ \delta^\al_1
\delta_{p,0}\,x
\eeq
for some formal power series $c^{\al,p}(\ve)$. We assume that only
finitely many of them are nonzero. The series must satisfy the condition of {\it
genericity} that we shall now formulate.

Let us expand the Hamiltonian densities (\ref{def-h}) in powers of $\ve$
\beq\label{theta}
h_{\al,p} = \theta_{\al,p+1}(v,u) + O(\ve).
\eeq
The explicit formulae for the functions $\theta_{\al,p}(v,u)$ can be found in
\cite{DZ3}. Let us impose the following assumptions for the leading terms of the
series $c^{\al,p}(\ve)$.

1. There exist values $\bar v$, $\bar u$ such that
\eqa\label{bar-uv}
&&
\sum_{p\geq 0} c^{\alpha,p}(0) \left.{\pal \theta_{\alpha,p}(v,u)\over \pal
v}\right|_{v=\bar v, u=\bar u}=0
\nn\\
&&
\sum_{p\geq 0} c^{\alpha,p}(0) \left.{\pal \theta_{\alpha,p}(v,u)\over \pal
u}\right|_{v=\bar v, u=\bar u}=0
\nn\\ 
\eeqa
and

2. The operator of multiplication by the vector
\beq\label{vektor}
\nabla \sum_{p\geq 1} c^{\al,p}(0) \theta_{\al,p-1}(\bar v, \bar u)
\eeq
is invertible element of the Frobenius algebra $T_{\bar v, \bar u}M_{\rm Toda}$.

Under these assumptions the following Lemma holds true (cf. \cite{DZ3}).

\begin{lemma} There exists a unique solution to the Euler - Lagrange equations
(\ref{EL-toda}) in the class of formal series
\eqa\label{formal}
&&
v=v(x, {\bf t}, \ve) = a_0(\ve) + \sum_{k>0} a_{\al_1, p_1; \dots ; \al_k,
p_k}(\ve) t^{\al_1, p_1} \dots t^{\al_k, p_k}|_{t^{1,0}\mapsto t^{1,0} + x}
\nn\\
&&
u=u(x, {\bf t}, \ve) = b_0(\ve) + \sum_{k>0} b_{\al_1, p_1; \dots ; \al_k,
p_k}(\ve) t^{\al_1, p_1} \dots t^{\al_k, p_k}|_{t^{1,0}\mapsto t^{1,0} + x}
\eeqa
where
\beq\label{formal0}
a_0(0)=\bar v, \quad b_0(0) = \bar u.
\eeq
\end{lemma}

\pf In the leading order in $\ve$ the Euler - Lagrange equations (\ref{EL-toda})
become just equations for the critical points of the function
$$
\sum_{p\geq 0} (t^{\al,p}-c^{\al,p}(0)+x\delta^\al_1 \delta^p_0) \theta_{\al,p}(v,u).
$$
The above two assumptions imply that there exists a unique critical
point
$$
v_0=v_0(x, {\bf t}, \ve), ~~u_0=u_0(x, {\bf t}, \ve)
$$
of the function
$$
\sum_{p\geq 0} (t^{\al,p}-c^{\al,p}(\ve)+x\delta^\al_1 \delta^p_0) \theta_{\al,p}(v,u)
$$
in the class of formal power series of the structure similar to (\ref{formal})
with
$$
v_0(0, {\bf 0}, 0) = \bar v, ~~ u_0(0, {\bf 0}, 0) = \bar u.
$$
It is easy to see that these functions can be uniquely extended to a 
solution to the full Euler - Lagrange equations
(\ref{EL-toda}). The Lemma is proved. \epf

\begin{lemma} The space of solutions of the Euler-Lagrange equation 
(\ref{EL-toda})
is invariant with respect to the flows of the extended Toda hierarchy.
\end{lemma}
\pf Let us represent the difference operators $A_{\beta,q}$ defined in (\ref{def-a-1}) 
by
\beq\label{v-a}
A_{\beta,q}=\sum_{k\ge 0} a_{\beta,q;k}\Lambda^k,\quad \beta=1,2,\ q\ge 0.
\eeq
Then from the definition of the Hamiltonians $H_{\beta,q}$
we obtain
\beq\label{toda-H4}
\frac{\delta H_{\al,p}}{\delta v (x)}=a_{\al,p;0}(x),\quad
\frac{\delta H_{\al,p}}{\delta u(x)}=a_{\al,p;1}(x-\ve)\,e^{u(x)}
\eeq
The Lax pair representation (\ref{td-flow-Lax}) of the extended Toda hierarchy yields
\eqa
&&\frac{\pal A_{\al,p}}{\pal t^{\beta,q}}-\frac{\pal A_{\beta,q}}{\pal t^{\al,p}} =[A_{\beta,q}, A_{\al,p}],\nn\\
&&\frac{\pal e^{u}}{\pal t^{\beta,q}}=\left[a_{\beta,q;0}(x)-a_{\beta,q;0}(x-\ve)\right] e^{u(x)}.\nn
\eeqa
These equations together with (\ref{v-a}) imply
\beq
\frac{\pal}{\pal t^{\beta,q}}\left(\frac{\delta H_{\al,p}}{\delta w^\gamma (x)}\right)
=\frac{\pal}{\pal t^{\al,p}}\left(\frac{\delta H_{\beta,q}}{\delta
w^\gamma(x)}\right), \quad \al,\beta,\gamma=1,2;\ p,q\ge 0.
\eeq
So, under the flows of the extended Toda hierarchy we have
\beq\label{toda-H5}
\frac{\pal}{\pal t^{\beta,q}}\left(\sum_{p\ge 0} {\tilde t}^{\al,p} \frac{\delta H_{\al,p-1}}
{\delta w^\gamma (x)}\right) =\frac{\delta H_{\beta,q-1}}{\delta w^\gamma (x)}+
\sum_{m\ge 0} \frac{\pal}{\pal w^{\xi,m}}\left(
\frac{\delta H_{\beta,q}}{\delta w^\gamma (x)} \right)
\pal_x^m \left(\sum_{p\ge 1} {\tilde t}^{\al,p} \frac{\pal w^\xi}{\pal t^{\al,p-1}} \right).
\eeq
Here and below we use the following notations for the $x$-derivatives of the
functions $u$ and $v$
$$
w^{\xi, m} := {\pal^m w^\xi\over \pal x^m}, \quad \xi = 1, \, 2, \quad m\geq 0.
$$
So
$$
w^{1,m} = v^{(m)}, \quad w^{2,m}= u^{(m)}.
$$
By using (\ref{string}) that we will prove in the Lemma \ref{lemma33} below we know that
the r.h.s. of (\ref{toda-H5}) can be rewritten as
$$
\frac{\delta H_{\beta,q-1}}{\delta w^\gamma (x)}-\frac{\pal}{\pal v}
\left(\frac{\delta H_{\beta,q}}{\delta w^\gamma (x)}\right)
$$
which equals zero due to (\ref{toda-H4}), (\ref{toda-H3}) and the identity $\frac{\pal L}{\pal v}=1$.
The Lemma is proved.
\epf

Due to the uniqueness of solutions of the initial value problem for the Euler-Lagrange
equation (\ref{EL-toda}) and the above theorem, we have
\begin{theorem}\label{td-el}
Any solution of the equation (\ref{EL-toda}) gives a solution to the extended Toda hierarchy.
\end{theorem}

Using quasitriviality it can be shown that the class of solutions of the extended Toda hierarchy that is given by the
above theorem form a dense subset of the class of its analytic solutions $w^\al(x,{\bf t},\ve),\, \al=1,2$
(see \cite{DZ3}). We call this class of solutions the generic class of solutions of the extended
Toda hierarchy, and we will restrict ourselves to it henceforth.
\begin{lemma}\label{lemma33}
Any solution $(v,u)$ of the Euler-Lagrange equation (\ref{EL-toda}) satisfies the equations
\eqa\label{string}
&&
\sum_{p\ge 1}{\tilde t}^{\al,p} \frac{\pal v}{\pal t^{\al,p-1}}+1=0
\nn\\
&&
\sum_{p\ge 1}{\tilde t}^{\al,p} \frac{\pal u}{\pal t^{\al,p-1}}=0
\\
&&
\sum_{q\ge 1} \left[ q\, \left({\tilde t}^{1,q}
\frac{\pal v}{\pal t^{1,q}}+{\tilde t}^{2,q-1}
\frac{\pal v}{\pal t^{2,q-1}}\right) + 2 \tilde t^{1,q} 
\frac{\pal v}{\pal t^{2,q-1}}\right] +v=0
\nn\\
&&
\sum_{q\ge 1} \left[ q\, \left({\tilde t}^{1,q}
\frac{\pal u}{\pal t^{1,q}}+{\tilde t}^{2,q-1}
\frac{\pal u}{\pal t^{2,q-1}}\right) + 2 \tilde t^{1,q} 
\frac{\pal u}{\pal t^{2,q-1}}\right] +2=0.
\label{Lm0s}
\eeqa
\end{lemma}
\pf The equation (\ref{string}) is the result of the application
of the operators $\frac{\pal}{\pal t^{2,0}}$ and $\frac1{\ve}
(\ld-1)$ on the Euler-Lagrange equation (\ref{EL-toda}). To prove the equation (\ref{Lm0s}),
we need to use the following bihamiltonian recursion relation of the extended Toda hierarchy \cite{CDZ}:
\beq\label{td-recur}
U_2^{\al\gamma}\,\frac{\delta H_{\beta,q-1}}{\delta w^\gamma}
=(q+\mu_\beta+\frac12)U_1^{\al\gamma}\,\frac{\delta H_{\beta,q}}{\delta w^\gamma}
+R^\gamma_\beta\, U_1^{\al\xi}\,\frac{\delta H_{\gamma,q-1}}{\delta w^\xi}.
\eeq
Here the first Hamiltonian structure $U^{\al\beta}_1$ is defined in (\ref{fpb})
and the second one is given by
\eqa\label{secpb}
&&U^{11}_2={1\over\epsilon}\left(e^{\epsilon\,\pal_x}e^{u(x)}-
e^{u(x)} e^{-\epsilon\,\pal_x}\right),\quad U^{12}_2={1\over \epsilon}
v(x)\left(e^{\epsilon\,\pal_x}-1 \right),\nn\\
&&U^{21}_2={1\over \epsilon}
\left(1-e^{-\epsilon\,\pal_x} \right)v(x),\quad
U^{22}_2={1\over \epsilon} \left(
e^{\epsilon\,\pal_x}-e^{-\epsilon\,\pal_x}\right).
\eeqa
The matrices $R$ and $\mu$ are defined by (\ref{def-mu-R}) (see also
(\ref{mu-r})).
Then the equation (\ref{Lm0s}) is obtained by applying the
operator $U_2^{\al\gamma}$ to both sides of the Euler-Lagrange
equation (\ref{EL-toda}) and by using the bihamiltonian recursion
relation (\ref{td-recur}).
The Lemma is proved. \epf

Let $v, u$ be any solution of the Euler-Lagrange
equation (\ref{EL-toda}) specified by a choice of the series $c^{\al,p}(\ve)$
and of the leading term $\bar v$, $\bar u$ in (\ref{bar-uv}). Due to Theorem \ref{thint} and theorem \ref{td-el} this solution
can be obtained from a solution $v_0, u_0$ of the dispersionless
Toda hierarchy. Denote by $\tau^{[0]}$ and $\tau$ the corresponding tau functions
with the relation
\beq\label{genus-b}
\log\tau=\ve^{-2}\,\log\tau^{[0]}+\sum_{g\ge 1} \ve^{2g-2}
F_g(w_0,\dots,w_0^{(3g-2)}).
\eeq
Note that the genus zero tau function $\tau^{[0]}$ is defined up to multiplication
by a function of the form $e^{\sum c^{[0]}_{\al,p} t^{\al,p}}$ with constants $c^{[0]}_{\al,p}$.
We now fix this ambiguity by taking
\beq\label{tau-00}
\log\tau^{[0]}=\frac12\,\sum \Omega^{[0]}_{\al,p;\beta,q}(v_0,u_0) {\tilde t}^{\al,p}\,{\tilde t}^{\beta,q}.
\eeq
where
\beq
\Omega^{[0]}_{\al,p;\beta,q}=\left.\Omega_{\al,p;\beta,q}\right|_{\ve=0}.
\eeq
The validity of this definition for the tau function of the solution $v_0, u_0$
of the dispersionless extended Toda hierarchy is based on the fact that $v_0, u_0$
satisfy the genus zero Euler-Lagrange equation
\beq
\sum {\tilde t}^{\al,p} \frac{\pal h_{\al,p-1}^{[0]}(v_0, u_0)}{\pal v_0}=0,
\quad 
\sum {\tilde t}^{\al,p} \frac{\pal h_{\al,p-1}^{[0]}(v_0, u_0)}{\pal u_0}=0
\eeq
with $h_{\al,p-1}^{[0]}=\left.h_{\al,p-1}\right|_{\ve=0}$. From this equation it readily follows that
\beq
\Omega^{[0]}_{\al,p;\beta,q}=\frac{\pal^2\log\tau^{[0]}}{\pal t^{\al,p}\pal t^{\beta,q}},
\eeq
and that the genus zero tau function satisfies the string equation
\beq\label{string-0}
\sum_{p\ge 1} {\tilde t}^{\al,p} \frac{\pal\log\tau^{[0]}}{\pal t^{\al,p-1}}+{\tilde t}^{1,0} {\tilde t}^{2,0}=0.
\eeq
The proof of the above statement can be found in \cite{D0}. It was proved in \cite {DZ2} that such a
tau function also satisfies the genus zero Virasoro constraints given 
by the Virasoro operators (\ref{def-vir-td}).
The action of these operators on tau functions of the form (\ref{genus-b}) can be expressed as
\beq
L_m(\ve^{-1}\tilde{\bf t}, \ve {\pal\over\pal{\bf t}})\tau=\left(\sum_{g\ge 0} \ve^{2g-2} Z_g\right) \tau,\quad m\ge -1.
\eeq
The genus zero Virasoro constraints are given by $Z_0=0$. We are to prove below
that the tau function
of a generic solution to the extended Toda hierarchy satisfies the full genera Virasoro
constraints $Z_g=0,\ g\ge 0$. Let us begin with the $L_{-1}$ and $L_0$
constraints.
\begin{lemma}
The tau function (\ref{genus-b}) satisfies the constraints
\beq\label{Lm1}
L_{-1}(\ve^{-1}\tilde{\bf t}, \ve {\pal\over\pal{\bf t}})\tau=0.\quad 
L_0(\ve^{-1}\tilde{\bf t}, \ve {\pal\over\pal{\bf t}})\tau=c_0(\ve)=\sum_{g\ge 1} \ve^{2g-2} c^{[g]}_0
\eeq
with certain constant $c_0(\ve)$.
\end{lemma}
\pf Let us apply the operator $\frac{\pal^2}{\pal t^{\sigma,k}\pal t^{\rho,l}}$
to the l.h.s. of the first equation of (\ref{Lm1}). Using the definition (\ref{def-omega})
for the tau function and the equation (\ref{string}) we get
\eqa
&&\ve^2\,\frac{\pal^2}{\pal t^{\sigma,k}\pal t^{\rho,l}}\left(\sum_{p\ge 1}{\tilde t}^{\al,p}
\frac{\pal \log\tau}{\pal t^{\al,p-1}}
+\frac1{\ve^2}\,{\tilde t}^{1,0} {\tilde t}^{2,0}\right)\nn\\
&&=\Omega_{\sigma,k-1;\rho,l}+\Omega_{\sigma,k;\rho,l-1}+
\sum_{p\ge 1}{\tilde t}^{\al,p}\frac{\pal\Omega_{\sigma,k;\rho,l}}
{\pal t^{\al,p-1}}+\left(\delta_{\sigma,1}\delta_{\rho,2}+
\delta_{\sigma,2}\delta_{\rho,1}\right)\delta_{k,0}\delta_{l,0} \nn\\
&&=\Omega_{\sigma,k-1;\rho,l}+\Omega_{\sigma,k;\rho,l-1}+
\sum_{p\ge 1}{\tilde t}^{\al,p}\frac{\pal\Omega_{\sigma,k;\rho,l}}
{\pal w^{\xi,m}}\pal_x^m \left(\frac{\pal w^\xi}{\pal
t^{\al,p-1}}\right)+ \left(\delta_{\sigma,1}\delta_{\rho,2}+
\delta_{\sigma,2}\delta_{\rho,1}\right)\delta_{k,0}\delta_{l,0}\nn\\
&&=\Omega_{\sigma,k-1;\rho,l}+\Omega_{\sigma,k;\rho,l-1}-
\frac{\pal\Omega_{\sigma,k;\rho,l}}{\pal v}
+\left(\delta_{\sigma,1}\delta_{\rho,2}+
\delta_{\sigma,2}\delta_{\rho,1}\right)\delta_{k,0}\delta_{l,0}\nn\\
&&=0.\label{eqtozero}
\eeqa
Here the last equality is due to (\ref{omega-du1}). On the other hand, the Euler-Lagrange equation
(\ref{EL-toda}) implies that the l.h.s. of the first formula of (\ref{Lm1}) does not
depend on ${\tilde t}^{1,0}$ and ${\tilde t}^{2,0}$, so there
exist constants $c_{-1}^{[g]},\ c_{\al,p}^{[g]}, \ \al=1,2,\ p\ge 0,\ g\ge 1$ such that
\beq
\sum_{p\ge 1}{\tilde t}^{\al,p} \frac{\pal \log\tau}{\pal
t^{\al,p-1}} +\frac1{\ve^2}\,{\tilde t}^{1,0} {\tilde
t}^{2,0}=\sum_{p\ge 1, g\ge 1} \ve^{2g-2}\,c^{[g]}_{\al,p-1} {\tilde
t}^{\al,p}+\sum_{g\ge 1}\ve^{2g-2}\, c^{[g]}_{-1}.\label{Lm1-b}
\eeq
Here the vanishing of the $\ve^{-2}$ term in the r.h.s. of the above identity
is due to (\ref{string-0}).
Thus if we modify the tau function by
\beq
\tau\mapsto {\tilde \tau}=\tau\,e^{\sum_{p\ge 0,g\ge 1} c^{[g]}_{\al,p}\,{\tilde t}^{\al,p}},
\eeq
then we obtain
\beq\label{st-1}
{L}_{-1}(\ve^{-1}\tilde{\bf t}, \ve {\pal\over\pal{\bf t}})\, {\tilde \tau}=c_{-1}(\ve)\,{\tilde\tau}=(\sum_{g\ge 1}\ve^{2g-2}\, c^{[g]}_{-1}){\tilde \tau}.
\eeq
We will prove the vanishing of the constants $c^{[g]}_{\al,p},\ c_{-1}(\ve)$ in a moment.

By using the formula (\ref{omega-dE1}) and a similar argument as that
given in the proof of (\ref{eqtozero}), we can prove the validity of the following identity
\beq
\frac{\pal^2}{\pal t^{\sigma,k}\pal t^{\rho,l}}\left(\frac{L_0{\tilde\tau}}{{\tilde \tau}}\right)=0.
\eeq
Here
$$
L_0 = L_0(\ve^{-1}\tilde{\bf t}, \ve {\pal\over\pal{\bf t}}).
$$
So there exist constants $c_0(\ve)$ and $b_{\al,p}(\ve)$ such that
\beq\label{w-const}
\frac{L_0{\tilde\tau}}{{\tilde \tau}}=\sum_{\al,p} b_{\al,p}(\ve) {\tilde t}^{\al,p}+c_0(\ve).
\eeq
By using the commutation relation $[L_{-1},L_{0}]=-L_{-1}$ we obtain
\beq\label{Lm10-b}
{L}_{-1} \left[\left(\sum b_{\al,p}(\ve) {\tilde t}^{\al,p}+c_0(\ve)\right){\tilde \tau}\right]-
L_0\left(c_{-1}(\ve)\,
{\tilde \tau}\right)=-c_{-1}(\ve) {\tilde\tau}.
\eeq
The l.h.s. of the above equality reads
\beq
\left(\sum_{p\ge 1} b_{\al,p-1} {\tilde t}^{\al,p}\right){\tilde \tau}
+\left(\sum b_{\al,p} {\tilde
t}^{\al,p}+c_0(\ve)\right) {L}_{-1}\,{\tilde \tau}-c_{-1}(\ve) {L}_0
{\tilde \tau} =\left(\sum_{p\ge 1} b_{\al,p-1} {\tilde t}^{\al,p}\right){\tilde\tau}.
\eeq
So from (\ref{Lm10-b}) it follows that
\beq
\left(\sum_{p\ge 1} b_{\al,p-1} {\tilde t}^{\al,p}\right){\tilde\tau}=-c_{-1}(\ve)\, {\tilde\tau},
\eeq
from which we obtain $c_{-1}(\ve)=b_{\al,p}(\ve)=0$.

Now we proceed to proving the vanishing of the constants $c^{[g]}_{\al,p}$. From the above argument we
already have the identity
$$
L_{-1}(\ve^{-1}\tilde{\bf t}, \ve {\pal\over\pal{\bf t}})\left(\tau\,e^{\sum_{g\ge 1} \epsilon^{2g} c^{[g]}_{\alpha,p} t^{\alpha,p}}\right)=0.
$$
At the genus one level we have
\beq\label{st-2}
\sum_{p\ge 1} {\tilde t}^{\al,p}\frac{\pal F_1(w,w_x)}{\pal t^{\al,p-1}}+
\sum_{p\ge 1} {\tilde t}^{\al,p}c^{[g]}_{\alpha,p-1}=0.
\eeq
Starting from this formula till the end of the proof of the Lemma we will
redenote for the sake of brevity the arguments $w_0=(v_0,u_0)$ and
${w_0}_x=({v_0}_x, {u_0}_x)$ of the function $F_1(w_0, {w_0}_x)$ by
$w=(v,u)$ and $w_x=(v_x, u_x)$.

Since $\tau^{[0]}$ satisfies the genus zero Virasoro constraints, we can use the vanishing of the
genus zero Virasoro symmetries to obtain, as we did in \cite{DZ2,DZ3}, the following formula
\beq
\sum_{p\ge 1} {\tilde t}^{\al,p}\frac{\pal F_1(w,w_x)}{\pal t^{\al,p-1}}=
-\frac{\pal F_1}{\pal v}.
\eeq
Thus the identity (\ref{st-2}) can be rewritten as
\beq
\sum_{p\ge 1} {\tilde t}^{\al,p}c^{[1]}_{\alpha,p-1}=\frac{\pal F_1}{\pal v}.
\eeq
By applying the operator $\sum_{p\ge 0} z^p\frac{\pal}{\pal t^{\al,p}}$ to the above identity we get
\eqa
&&\sum_{p\ge 1} c^{[1]}_{\al,p-1} z^p=\sum_{p\ge 0}\left[\frac{\pal^2 F_1}{\pal v\pal w^\gamma}
\pal_x\left(\frac{\pal\theta_{\al,p+1}}{\pal w_\gamma}\right)+
\frac{\pal^2 F_1}{\pal v\pal w^\gamma_x}
\pal_x^2\left(\frac{\pal\theta_{\al,p+1}}{\pal w_\gamma}\right)\right]\,z^p\nn\\
&&=\left[\frac{\pal^2 F_1}{\pal v\pal w^\gamma} c^{\gamma\sigma}_\rho w^\rho_x+
\frac{\pal^2 F_1}{\pal v\pal w^\gamma_x}\,\pal_x(c^{\gamma\sigma}_\rho w^\rho_x)+
z\,\frac{\pal^2 F_1}{\pal v\pal w^\gamma_x}\,c^{\gamma\sigma_1}_{\rho_1} w^{\rho_1}_x\,
c^\sigma_{\sigma_1\rho_2} w^{\rho_2}_x\right]\,\frac{\pal\theta_\al(z)}
{\pal w^\sigma}.
\nn\\
\label{d-s}
\eeqa
Here the functions $\theta_{\al,p}=\theta_{\al,p}(w), \
c_{\al\beta\gamma}=c_{\al\beta\gamma}(w)$ are given by
\eqa
&&\theta_{\al,p}=h_{\al,p-1}|_{\ve=0},\quad \theta_\al(z)=\sum_{p\ge 0} \theta_{\al,p} z^p,\nn\\
&&c_{\al\beta\gamma}=\frac{\pal^3}{\pal w^\al\pal w^\beta\pal w^\gamma}
(\frac12\,v^2 u+e^{u}).
\eeqa
and the raise of indices in $c_{\al\beta\gamma}$ is done by the metric
(\ref{metric}), i.e.
$\eta^{11}=\eta^{22}=0, \eta^{12}=\eta^{21}=1$.
In the above computation we used the horizontality of the differentials
of the functions $\theta_\alpha(w;z)$ w.r.t. the deformed flat connection on
$M_{\rm Toda}$, i.e. the equations
\beq
\frac{\pal^2\theta_\al(z)}{\pal w^\beta\pal w^\gamma}=
z c^\xi_{\beta\gamma}\,\frac{\pal \theta_{\al}(z)}{\pal w^\xi}.
\eeq
So from (\ref{d-s}) we get
\eqa
&&\frac{\pal^2 F_1}{\pal v\pal w^\gamma} c^{\gamma\sigma}_\rho w^\rho_x+
\frac{\pal^2 F_1}{\pal v\pal w^\gamma_x}\,\pal_x(c^{\gamma\sigma}_\rho w^\rho_x)=0,\nn\\
&&\left(\frac{\pal^2 F_1}{\pal v\pal w^\gamma_x}\,c^{\gamma\sigma_1}_{\rho_1} w^{\rho_1}_x\,
c^\sigma_{\sigma_1\rho_2} w^{\rho_2}_x\right) \eta_{\sigma\al}=c^{[1]}_{\al,0}.
\eeqa
and together with (\ref{d-s}) these formulae in turn yield
\beq\label{1-a}
\sum_{p\ge 0} c_{\al,p}^{[1]} z^{p}=c^{[1]}_{\gamma,0}\,\pal^\gamma\theta_{\al}(z)
\eeq
By differentiating both sides of the above equation w.r.t. $x$ we get
\beq
0=c_{\gamma,0}^{[1]}\,c^{\gamma}_{\xi\sigma}\,v^\sigma_x\,\pal^\xi\theta_\al(z)
\eeq
which implies
$$
c_{\gamma,0}^{[1]}=0.
$$
So from (\ref{1-a}) we obtain
$$
c_{\gamma,p}^{[1]}=0, \quad p\ge 1.
$$

In a completely similar way we can prove that
$$
c_{\gamma,p}^{[g]}=0, \quad p\ge 0, \ g\ge 2.
$$
The Lemma is proved. \epf

The following Lemma represents the bihamiltonian recursion relation for the extended
Toda hierarchy in terms of its tau functions:

\begin{lemma}\label{lemma-1}
The following recursion relations hold true for any $q\geq 1$
for the tau functions of generic solutions to the extended Toda hierarchy:
\eqa\label{recur-tau1}
&&
q\, \left(\Lambda-1\right) {\pal \log\tau\over \pal t^{1,q}}
={\cal R} {\pal \log\tau\over \pal t^{1,q-1}}
- 2 \left(\Lambda-1\right) {\pal \log\tau\over\pal t^{2,q-1}}
\\
&&
(q+1)\, \left(\Lambda-1\right) {\pal \log\tau\over \pal t^{2,q}}
={\cal R} {\pal \log\tau\over \pal t^{2,q-1}}\label{recur-tau2}
\eeqa
where the operator ${\cal R}$ is defined by 
\beq\label{def-recur-R}
{\cal R}= v(x) (\Lambda-1) +\epsilon
(\Lambda+1) {\partial\over\partial t^{2,0}}.
\eeq
\end{lemma}
\pf Denote
\beq
W_{\beta,q}:={\cal R} {{\partial \log \tau}\over{\partial t^{\beta,q-1}}}- (q+\mu_\beta+\frac12)
(\Lambda-1) {{\partial \log \tau}\over{\partial
t^{\beta,q}}} -(\Lambda-1) R^\gamma_\beta
{{\partial \log \tau} \over{\partial t^{\gamma,q-1}}}.
\eeq
We are to prove that $W_{\beta,q}=0$ for $\beta=1,\, 2$, $q\geq 1$.
From Lemma \ref{prop-1} and from the bihamiltonian recursion relation (\ref{td-recur})
with $\al=2$ we obtain by a direct calculation that
\beq\label{foll}
(\Lambda-1) W_{\beta,q}=0.
\eeq
We note that $W_{\beta,q}$ can be expressed as 
homogeneous differential polynomials
in $w^{i,m},\newline e^{\pm u}, i=1,2, m\ge 0$  of 
degree $q+\frac32+\mu_\beta$. Recall that the degree of
such differential polynomials is defined in (\ref{grad-A}). 
So the Lemma follows from the above equation (\ref{foll}).
The Lemma is proved.\epf

\noindent{\bf Proof of the Main Theorem}
Let us first prove that the following recursion relation holds true:
\beq\label{recur-vir}
{\cal R} \frac{L_m \tau}{\tau}=(\ld-1)\frac{L_{m+1} \tau}{\tau}.
\eeq
Here and below
$$
L_m = L_m(\ve^{-1}\tilde {\bf t}, \ve{\pal\over\pal {\bf t}}).
$$
From the definition of the operator ${\cal R}$ we have
\eqa &&{\cal R}\sum_{k=1}^{m-1} k! (m-k)!\left(\frac1{\tau}\,\frac{\pal^2\tau}{\pal t^{2,k-1} \pal
t^{2,m-1-k}}\right) \nn\\
&&= \sum_{k=1}^{m-1} k!(m-k)!\left[ \frac{\pal}{\pal t^{2,k-1}} {\cal
R}\frac{\pal\log\tau}{\pal t^{2,m-1-k}}+\left({\cal
R}\frac{\pal\log\tau}{\pal
t^{2,k-1}}\right)\,\ld\,\frac{\pal\log\tau}{\pal
t^{2,m-1-k}}\right.\nn\\
&&\left. +\frac{\pal\log\tau}{\pal t^{2,k-1}}\,{\cal
R}\frac{\pal\log\tau}{\pal t^{2,m-1-k}}\right],\quad m\ge -1.\nn
\eeqa
So by using the recursion relations (\ref{recur-tau1}), (\ref{recur-tau2}) 
we can deduce
the relation (\ref{recur-vir}) for $m\ge 1$ as follows:
\eqa &&{\cal R} \frac{{L_m}\tau}{\tau}= \ve^2 (\ld-1)
\sum_{k=1}^{m} k! (m+1-k)!\,\frac1{\tau}\,\frac{\pal^2\tau}{\pal
t^{2,k-1} \pal
t^{2,m-k}} \nn\\
&& -\ve^2 m! (\ld-1) \frac{\pal^2\log\tau}{\pal t^{2,m-1} \pal
t^{2,0}} -m!\left[\ve (\ld-1)\frac{\pal\log\tau}{\pal
t^{2,0}}\right]
\ve (\ld+1)\frac{\pal\log\tau}{\pal t^{2,m-1}}\nn\\
&&+\sum_{k\ge 1} \frac{(m+k)!}{(k-1)!} (\ld-1) \left[(m+k+1)\left(
{\tilde t}^{1,k}\frac{\pal\log\tau}{\pal t^{1,m+k+1}}+ {\tilde
t}^{2,k-1}\frac{\pal\log\tau}{\pal t^{2,m+k}}\right)\right.
\nn\\
&&\quad\left. +2{\tilde t}^{1,k} \frac{\pal\log\tau}{\pal
t^{2,m+k}}
\right]+\ve (m+1)! (\ld+1) \frac{\pal\log\tau}{\pal t^{2,m}}\nn\\
&&+2\sum_{k\ge 0} (m+k+1) \al_m(k) {\tilde t}^{1,k} (\ld-1)
\frac{\pal\log\tau}{\pal t^{2,m+k}}\nn\\
&&+2 \ve \al_m(0) v\ld  \frac{\pal\log\tau}{\pal t^{2,m-1}}
+2 \ve^2 \al_m(0) \ld \frac{\pal^2\log\tau}{\pal t^{2,m-1}\pal t^{2,0}}\nn\\
&&=(\ld-1) \frac{{L_{m+1}}\tau}{\tau} -\ve^2 m! (\ld-1)
\frac{\pal^2\log\tau}{\pal t^{2,m-1} \pal
t^{2,0}}\nn\\
&& -m!\left[\ve (\ld-1)\frac{\pal\log\tau}{\pal t^{2,0}}\right]
\ve (\ld+1)\frac{\pal\log\tau}{\pal t^{2,m-1}}\nn\\
&&+\ve (m+1)! (\ld+1) \frac{\pal\log\tau}{\pal t^{2,m}}
-2 \ve (m+1)!\ld  \frac{\pal\log\tau}{\pal t^{2,m}}\nn\\
&&+2 \ve \al_m(0) v\ld  \frac{\pal\log\tau}{\pal t^{2,m-1}}
+2 \ve \al_m(0) \ld \frac{\pal^2\log\tau}{\pal t^{2,m-1}\pal t^{2,0}}\nn\\
&&=(\ld-1) \frac{{L_{m+1}}\tau}{\tau} +\ve m!\, {\cal R}
\frac{\pal\log\tau}{\pal t^{2,m-1}}
-\ve (m+1)! (\ld-1) \frac{\pal\log\tau}{\pal t^{2,m}}\nn\\
&&=(\ld-1) \frac{{L_{m+1}}\tau}{\tau}. \eeqa It can be
easily checked that the recursion relation (\ref{recur-vir}) is
also true for $m=-1, 0$.

From (\ref{Lm1}) and the recursion relation (\ref{recur-vir}) we know that
\beq\label{d-s-2}
(\ld-1)\left(\frac{L_{1} \tau}{\tau}\right)=0.
\eeq
Since the function $\tau^{[0]}$ satisfies the genus zero Virasoro constraints, it follows from
(\ref{genus-b}) that
\beq
\frac{L_1\tau}{\tau}=\sum_{g\ge 1}\ve^{2g-2} W_g(w_0,w_0 ',\dots,w_0^{(m_g)}).
\eeq
Thus from (\ref{d-s-2}) we arrive at the formula
\beq\label{d-s-3}
\frac{L_1\tau}{\tau}=c_1(\ve)=\sum_{g\ge 1}\ve^{2g-2} c_1^{[g]}
\eeq
for some constants $c^{[g]}_1$.

On the other hand, by using the  the commutation relation (\ref{comm-L}) and the
equations (\ref{Lm1}) we obtain
\beq\label{L-tau-act}
L_{-1}\left(L_m \tau\right)=-(m+1) L_{m-1}\tau, \quad
L_{0}\left(L_m \tau\right)=(c_0(\ve)-m) L_{m}\tau.
\eeq
These formulae can be rewritten as
\beq\label{452}
{\hat L}_{-1}\left(\frac{L_m \tau}{\tau}\right)=-(m+1)
\left(\frac{L_{m-1} \tau}{\tau}\right),\quad {\hat
L}_{0}\left(\frac{L_m \tau}{\tau}\right)=-m
\left(\frac{L_{m} \tau}{\tau}\right).
\eeq
Here ${\hat L}_{-1}=L_{-1}-\frac1{\ve^2} {\tilde t}^{1,0}\,{\tilde t}^{2,0}$ and
${\hat L}_0=L_{0}-\frac1{\ve^2} \left({\tilde t}^{1,0}\right)^2$. By
putting $m=1$ into the above two relations we get
\beq\label{c00}
{\hat L}_{-1}\left(\frac{L_1 \tau}{\tau}\right)=-2\, c_0(\ve),
\quad {\hat L}_{0}\left(\frac{L_1 \tau}{\tau}\right)=-
\left(\frac{L_{1} \tau}{\tau}\right).
\eeq
So from (\ref{d-s-3}) we have $c_0(\ve)=c_1(\ve)=0$, and we proved the vanishing of
$L_0\tau$ and $L_1\tau$.

By using the recursion relation (\ref{recur-vir}) with $m=1$ we obtain
\beq\label{d-s-5}
(\ld-1)\left(\frac{L_{2} \tau}{\tau}\right)=0.
\eeq
Due to the same reason as the one we used to derive (\ref{d-s-3}) we have
\beq\label{d-s-4}
\frac{L_2\tau}{\tau}=c_2(\ve)=\sum_{g\ge 1}\ve^{2g-2} c_2^{[g]}
\eeq
for some constants $c^{[g]}_2$.
So by using the second formula in (\ref{452}) we get $L_2\tau=0$. Now, the
Virasoro commutation relation (\ref{comm-L}) implies the validity of all the Virasoro
constraints 
$$
L_m(\ve^{-1}\tilde{\bf t}, \ve{\pal\over\pal{\bf t}})\tau=0,\ m\ge -1.
$$ 

It remains to prove that linear action of the Virasoro operators
(\ref{def-vir-td}) defines infinitesimal symmetries of the extended Toda
hierarchy. To this end we observe that
$$
L_m(\ve^{-1} ({\bf t}-{\bf c}(\ve)), \ve {\pal\over\pal{\bf t}})=
L_m (\ve^{-1} {\bf t}, \ve {\pal\over\pal{\bf t}})- a^{\al,p}{\pal\over\pal
t^{\al,p}} - b_{\al,p}t^{\al,p} -c
$$
where $a^{\al,p}$, $b_{\al,p}$ and $c$ are some series in $\ve$ that may also
depend on $m$. Note that the $a$ series contains only nonnegative powers
of $\ve$. From the already proven Virasoro constraints it follows that, for
any generic solution to the extended Toda hierarchy the action of the Virasoro
operators on the tau function can be recast as
$$
L_m (\ve^{-1} {\bf t}, \ve {\pal\over\pal{\bf t}})\tau = 
(a^{\al,p}{\pal\over\pal
t^{\al,p}} + b_{\al,p}t^{\al,p} +c)\tau.
$$
So, for a small parameter $\delta$
\beq\label{konec}
\tau + \delta\, L_m(\ve^{-1} {\bf t}, \ve {\pal\over\pal{\bf
t}})\tau
= e^{\delta (b_{\al,p}t^{\al,p} +c)} e^{\delta \,a^{\al,p}{\pal\over\pal
t^{\al,p}} } \tau +O(\delta^2).
\eeq
The operator 
$$
e^{\delta \,a^{\al,p}{\pal\over\pal
t^{\al,p}} }
$$
is nothing but the shift along trajectories of the hierarchy. Note that such a
shift leaves invariant the class of generic solutions. Multiplication by
the exponential of a linear function in the times for obvious reasons maps a tau
function to another one for the same solution to the hierarchy. 
This proves that (\ref{konec}) is again a tau function of some solution of the
extended Toda hierarchy.
The Theorem is proved.\epf

\setcounter{equation}{0} \setcounter{theorem}{0}
\section{The topological solution of the extended Toda hierarchy}

Let us briefly recall the definition of Gromov - Witten invaraints and their
descendents and the construction of the Witten's generating function (physicists
call it {\it free energy} of the two-dimensional $CP^1$ topological sigma model).
Denote 
$\phi_1=1\in H^0(CP^1)$, $\phi_2=\omega\in H^2(CP^1)$ the basis in the
cohomology space $H^*(CP^1)$.
The 2-form $\omega$ is assumed to be normalized by the condition
$$
\int_{CP^1}\omega=1.
$$

The free energy of the $CP^1$ topological sigma-model is a function of infinite number of {\it coupling parameters}
$$
{\bf t}=(t^{1,0}, t^{2,0}, t^{1,1}, t^{2,1}, \dots)
$$
and of $\ve$
defined by
the following genus expansion form:
\beq\label{free}
{\cal F}({\bf t}; \ve)=\sum_{g\ge 0}\ve^{2g-2} {\cal F}_g({\bf t}).
\eeq
The parameter $\ve$ is called here the string coupling constant, and the function ${\cal F}_g={\cal F}_g({\bf t})$
is called the genus $g$ free energy which is given by
\beq
{\cal F}_g=\sum_{} \frac1{m!} t^{\al_1,p_1}\dots t^{\al_m,p_m} \langle\tau_{p_1}(\phi_{\al_1})\dots
\tau_{p_m}(\phi_{\al_m})\rangle_g,
\eeq
where $\tau_p(\phi_\al)$ are the gravitational descendent of the primary fields 
$\phi_\al$, $t^{\al,p}$ is the corresponding coupling constants ,
and the rational numbers 
$\langle\tau_{p_1}(\phi_{\al_1})\dots \tau_{p_m}(\phi_{\al_m})\rangle_g$ 
are given by the following intersection numbers on the moduli spaces
of $CP^1$-valued stable curves of genus $g$:
\beq\label{korrel}
\langle\tau_{p_1}(\phi_{\al_1})\dots \tau_{p_m}(\phi_{\al_m})\rangle_g = \sum_{\beta} q^\beta\int_{[\bar M_{g,m}(CP^1,\beta)]^{\rm virt}}
{\rm ev}_1^* \phi_{\al_1} \wedge \psi_1^{p_1} \wedge \dots \wedge {\rm ev}_m^* \phi_{\al_m} \wedge \psi_m^{p_m}.
\eeq
Here $\bar M_{g,m}(CP^1,\beta)$ is the moduli space of stable curves of genus $g$ with $m$ markings of the given degree $\beta\in H_2 (CP^1;
{\mathbb Z})$, ${\rm ev}_i$ is the evaluation map
$$
{\rm ev}_i: \bar M_{g,m}(CP^1,\beta)\to CP^1
$$
corresponding to the $i$-th marking,
$\psi_i$ is the first Chern class of the tautological line bundle over the moduli space corresponding to the $i$-th marking. 
According to the divisor
axiom \cite{KM} the dependence of the Gromov - Witten potential
on the indeterminate $q$ 
appears only through the combination $q \, e^{t^{2,0}}$. We will therefore omit
the dependence on $q$ in the formulae.

Let us clarify the relationship between our theory of Virasoro symmetries
of the extended Toda hierarchy and the Virasoro conjecture of T.Eguchi, K.Hori,
and C.-S.Xiong \cite{eguchi,eguchi3} extended by S.Katz. Denote
$$
Z_{CP^1}({\bf t};\ve):= e^{{\cal F}({\bf t}; \ve)}
$$
the partition function of the $CP^1$ topological sigma-model. Here
${\cal F}({\bf t}; \ve)$ is the generating function 
of the $CP^1$ Gromov - Witten invariants and their descendents defined
above. According to the results of A.Givental \cite{givental1, givental2}
this partition function satisfies the following infinite sequence of linear
Virasoro constraints 
\eqa\label{giv-vir}
&&
L_{-1} Z_{CP^1}= {\pal \over\pal t^{1,0}}Z_{CP^1}
\nn\\
&&
L_m Z_{CP^1} = (m+1)! \left[ {\pal\over\pal t^{1,m+1}} + 2 \kappa_m
{\pal\over\pal t^{2,m}}\right] Z_{CP^1}, \quad m\geq 0
\eeqa
where
$$
\kappa_m=\sum_{j=1}^{m+1} \frac1j.
$$
Here $L_m$ are just the Virasoro operators defined in (\ref{def-vir-td}).
For the particular case of $CP^1$ (\ref{giv-vir}) coincide with the Virasoro
constraints conjectured in \cite{eguchi}. However, in their papers Eguchi, Hori
and Xiong formulated a somewhat more bold conjecture that says that all $g\geq
1$ Gromov - Witten invariants and their descendents of a smooth projective
variety can be {\it uniquely determined} by solving recursively the linear system of Virasoro
constraints. Although this conjecture seems to be too nice to be true in general
(Calabi - Yau
manifolds give counterexamples to uniqueness, see \cite{katz}), in certain cases
it can be justified.

Let us give our version of the Eguchi - Hori - Xiong Virasoro constraints
programme adapted to computing Gromov - Witten invariants of $CP^1$.

{\bf Step 1}. Computation of the genus zero Gromov - Witten potential
${\cal F}_0({\bf t})$. This can be done in terms of the Frobenius manifold
$M_{\rm Toda}$ as in \cite{npb, D1}. For the reader's convenience we
recall the algorithm of computation of ${\cal F}_0({\bf t})$ in the Appendix
below. Introduce functions
\eqa
&&
v_0 =v_0({\bf t}):= {\pal^2 {\cal F}_0({\bf t})\over \pal t^{1,0} \pal t^{2,0}}
\nn\\
&&
u_0=u_0({\bf t}):= {\pal^2 {\cal F}_0({\bf t})\over \pal t^{1,0} \pal t^{1,0}}
\nn
\eeqa
We will denote $u_0'$, $v_0'$, $u_0''$, $v_0''$ etc. the derivatives of these
functions along $t^{1,0}$.

{\bf Step 2}. Eguchi - Hori - Xiong $(3g-2)$-ansatz 
for the higher genus corrections.
Look for the $g\geq 1$ terms in the genus expansion (\ref{free}) in the form
\beq\label{ansatz}
{\cal F}_g({\bf t}) = F_g(v_0({\bf t}), u_0({\bf t}), v_0'({\bf t}), 
u_0'({\bf t}),
\dots, v_0^{(3g-2)}({\bf
t}),
u_0^{(3g-2)}({\bf t})), \quad g\geq 1.
\eeq
The ansatz (\ref{ansatz}) was proved by E.Getzler in \cite{getzler1}. In the
setup of our theory \cite{DZ3} of integrable systems the $(3g-2)$-ansatz
is a consequence of a deep result about quasitriviality of tau-symmetric
deformations of Poisson pencils.

{\bf Step 3}. Virasoro Conjecture. 

\noindent Part 1. The series (\ref{free}) where ${\cal F}_0({\bf t})$ is
the genus zero Gromov - Witten potential and the terms of positive genera
have the form (\ref{ansatz}) satisfies the Virasoro constraints (\ref{giv-vir}).
(Clearly the $(3g-2)$-ansatz is of no importance so far.)

\noindent Part 2. The degree $2g-2$ homogeneous
functions $F_g$  on the $(3g-2)$ jet space of $M_{\rm Toda}$ 
for all $g\geq 1$
are uniquely determined from the Virasoro constraints (\ref{giv-vir})
by solving recursively systems of linear equations.

Part 1 of the Virasoro Conjecture was proved by A.Givental \cite{givental1,
givental2}. Part 2 was proved in much more general framework of an arbitrary
semisimple Frobenius manifold in \cite{DZ3}. Combining these results we arrive
at

\begin{theorem} 1. The partition function $Z_{CP^1}({\bf t}; \ve)$ of the $CP^1$
topological sigma-model is uniquely determined by the Virasoro Conjecture
equations.

\noindent 2. It coincides with the tau function $\tau_{CP^1}$ of a particular solution
to the extended Toda hierarchy (\ref{td-flow-Lax})
\beq\label{same}
Z_{CP^1}({\bf t};\ve) =\tau_{CP^1}({\bf t};\ve)
\eeq
specified by the following choice of the
shift parameters $c^{\al,p}(\ve)$ and the initial point $\bar v$, $\bar u$:
\beq\label{topo}
c^{\al,p}(\ve) = \delta^\al_1 \delta^p_1, \quad \bar v = \bar u = 0.
\eeq
\end{theorem}

The choice (\ref{topo}) selects the solution satisfying the string equation
\beq
\sum_{p\ge 1} t^{\al,p}\frac{\pal {\cal F}}{\pal t^{\al,p-1}}+\frac1{\ve^2} t^{1,0}\,t^{2,0}=\frac{\pal{\cal F}}{\pal t^{1,0}}.
\eeq

{\bf Sketch of the proof}. As it was shown in \cite{DZ3}, from
validity of the Virasoro
constraints for the sum $\Delta{\cal F}$
of all $g\geq 1$ corrections to the Gromov - Witten
potential represented via $(3g-2)$-ansatz 
$$
\Delta {\cal F}: = \sum_{g\geq 1}\ve^{2\, g} F_g(v,u; v_x, u_x, \dots, v^{(3g-2)},
u^{(3g-2)})
$$
it follows the following {\it loop equation}
\eqa
&&
\sum_{r\geq 0}\left({\pal \Delta{\cal F}\over \pal v^{(r)}} \pal_x^r
{v-\lambda\over D} - 2 {\pal \Delta{\cal F}\over \pal u^{(r)}} 
\pal_x^r {1\over D}\right)
\nn\\
&&
+\sum_{r\geq 1} \sum_{k=1}^r \left(\matrix{r\cr k}\right) \pal_x^{k-1}
{1\over \sqrt{D}}\left( {\pal \Delta{\cal F}\over \pal v^{(r)}}
\pal_x^{r-k+1} {v-\lambda\over \sqrt{D}} -2 
{\pal \Delta{\cal F}\over \pal u^{(r)}}\pal_x^{r-k+1} {1\over \sqrt{D}}\right)
\nn\\
&&
=D^{-3} e^{u} \left( 4\, e^{u}+(v-\lambda)^2\right)
\nn\\
&&
+\sum_{k,l}{\epsilon^2\over 4} \left[ -\left( 
{\pal^2 \Delta{\cal F}\over \pal v^{(k)}\pal v^{(l)}}+
{\pal \Delta{\cal F}\over \pal v^{(k)}}
{\pal \Delta{\cal F}\over \pal v^{(l)}}\right)
\, \pal_x^{k+1} {v-\lambda\over \sqrt{D}} \pal_x^{l+1} {v-\lambda\over
\sqrt{D}}\right.
\nn\\
&&
+4\, \left( 
{\pal^2 \Delta{\cal F}\over \pal v^{(k)}\pal u^{(l)}}+
{\pal \Delta{\cal F}\over \pal v^{(k)}}
{\pal \Delta{\cal F}\over \pal u^{(l)}}\right)\,\pal_x^{k+1} {v-\lambda\over
\sqrt{D}} \pal_x^{l+1} {1\over \sqrt{D}}
\nn\\
&&
\left. -4\, \left( 
{\pal^2 \Delta{\cal F}\over \pal u^{(k)}\pal u^{(l)}}+
{\pal \Delta{\cal F}\over \pal u^{(k)}}
{\pal \Delta{\cal F}\over \pal u^{(l)}}\right)\,
\pal_x^{k+1} {1\over \sqrt{D}} \pal_x^{l+1} {1\over \sqrt{D}}\right]
\nn\\
&&
-{\epsilon^2\over 2}\sum_k \left\{ {\pal \Delta{\cal F}\over \pal v^{(k)}}
\pal_x^{k+1} e^{u}{4 \, e^{u}
(v-\lambda) u' -[(v-\lambda)^2 + 4\, e^{u}]\, v' \over D^3}
\right.
\nn\\
&&
\left.  +{\pal \Delta{\cal F}\over \pal u^{(k)}} \pal_x^{k+1}
e^{u} {4 \, (v-\lambda) \, v' - [ (v-\lambda)^2 + 4\, e^{u}]\, u'
\over D^3}\right\}
\nn\\
\eeqa
where
$$
D= (v-\lambda)^2 -4\, e^{u}.
$$
Here $\lambda$ is an arbitrary complex parameter. Expanding the loop equation
near $\lambda=\infty$ reproduces the Virasoro constraints for $\Delta
{\cal F}$. The proof of existence and uniqueness of the solution to this
equation is based on expanding the loop equation near zeroes 
$$
u_\pm =v \pm 2 e^{u/2}
$$
of $D$ (these are
the {\it canonical coordinates} on the Frobenius manifold $M_{\rm Toda}$).
The uniqueness of solution to the loop equation proves first part of the
Theorem.

To prove the second part we use the following arguments. From 
\cite{D1,Egu5, Egu4, Egu6}
we already know that 
$$
\tau^{[0]}({\bf t}):= {\cal F}_0({\bf t})
$$
is the tau function of the dispersionless extended Toda hierarchy. This solution
is specified by the shift parameters and the leading term (\ref{topo}).

The
transformation 
\beq\label{trans}
\log\tau^{[0]}\mapsto \log\tau^{[0]}+\Delta {\cal F}=: \ve^2\log\tau
\eeq
maps dispersionless tau functions to tau functions of the {\it full hierarchy}
associated with the semisimple Frobenius manifold $M_{\rm Toda}$. The full
hierarchy is uniquely determined, for the given semisimple Frobenius manifold
by the following properties:

\noindent - bihamiltonian structure satisfying certain nondegeneracy conditions;

\noindent - tau symmetry that provides existence of a tau function for a generic
solution; 

\noindent - invariance with respect to the linear action of the Virasoro
operators $L_m$, $m\geq -1$ onto the tau functions.

As we explained in the Introduction, the first two properties are met by the
extended Toda hierarchy due to results of \cite{CDZ}. The last property of
Virasoro invariance is established in the present paper. This implies that
the full hierarchy associated with the Frobenius manifold $M_{\rm Toda}$
coincides with the extended Toda hierarchy. Therefore the transformation
(\ref{trans}) maps the tau function $\tau^{[0]}$
of an arbitrary solution fo the dispersionless hierarchy
to the tau function $\tau$ of a solution of the full extended Toda hierarchy. Taking
$\tau^{[0]}={\cal F}_0$ one obtains $\tau=Z_{CP^1}$.
The Theorem is proved.\epf

Clearly the Theorem covers Corollaries \ref{cor-01} and \ref{cor-00}
formulated in the Introduction.

To illustrate the algorithm of computation of the genus expansion
(\ref{free}) for $CP^1$ let us write it down the first two terms of the
expansion. The formulae become simpler when written in the canonical coordinates
$$
u_\pm = v_0 \pm 2 e^{u_0\over 2}.
$$
Genus 1:
$$
F_1 =\frac1{24} \log u_+' u_-' -{1\over 12} \log {u_+ - u_-\over 4}.
$$
Genus 2:
\eqa
&&24^2\, F_2=
\frac{{4\,u_+''}^3\,(u_+-u_-)}{5\,{u_+'}^4} - 
  \frac{{4\,u_-''}^3\,(u_+-u_-)}{5\,{u_-'}^4}-
\frac{u_+''\,u_-''}{4\,u_+'\,u_-'}
\nn\\
&&\quad+
\frac{3\,u_+''}{4\,{u_+'}^3}\left[\frac12\,u_+''\,u_-'-\frac75\,u_+'''\,(u_+-u_-)
\right]
+
\frac{3\,u_-''}{4\,{u_-'}^3}\left[\frac12\,u_-''\,u_+'+\frac75\,u_-'''\,(u_+-u_-)
\right]\nn\\
&&\quad+
\frac1{4\,{u_+'}^2}\left[\frac{33}{10}\,{u_+''}^2-\frac9{10}\,u_+'''\,u_-'+
\frac1{10}\,u_+''\,u_-''+u_+^{IV}\,(u_+-u_-)\right]\nn\\
&&\quad+
\frac1{4\,{u_-'}^2}\left[\frac{33}{10}\,{u_-''}^2-\frac9{10}\,u_-'''\,u_+'+
\frac1{10}\,u_+''\,u_-''-u_-^{IV}\,(u_+-u_-)\right]\nn\\
&&\quad-
\frac1{4\,u_+'}\left(\frac{17}{5}\,u_+'''+\frac1{2}\,u_-'''\right)-
\frac1{4\,u_-'}\left(\frac{17}{5}\,u_-'''+\frac1{2}\,u_+'''\right)\nn\\
&&\quad-
\frac1{10\, (u_+-u_-)^2}\left(\frac{{u_+'}^3}{u_-'}+\frac{{u_-'}^3}{u_+'}
\right)
-\frac{1}{(u_+-u_-)^2}\left({u_+'}^2-\frac{11}{5}\,u_+'\,u_-'+{u_-'}^2\right)
\nn\\
&&\quad+
\frac{u_+''-u_-''}{u_+-u_-}\left(\frac{u_-'}{5\,u_+'}+\frac{u_+'}{5\,u_-'}+
1\right).
\eeqa

{\noindent \it Remark 1}. In \cite{getzler}, Getzler proved that, under 
the assumption of the recursion
relation (\ref{recur-tau2}), validity of the Virasoro constraints for 
$\tau_{CP^1}$ is equivalent
to (\ref{recur-tau1}). In his proof 
a recursion relation of the
form (\ref{recur-vir}) was used. The recursion (\ref{recur-tau2}) for 
$\tau_{CP^1}$ was proved in
\cite{oko2} on the subspace $\{t^{1,k}=0,\ k>1\}$ of the large phase space
of all couplings. Using this result
Getzler also proved (\ref{recur-tau1}) and (\ref{recur-tau2})
under the assumption of the Virasoro
constraints for $\tau_{CP^1}$. He did not consider connections between
recursion relations and Virasoro constraints for other solutions to the extended
Toda hierarchy.
Our Corollary \ref{cor-00} shows that the recursion
relations (\ref{recur-tau1}), (\ref{recur-tau2}) for $\tau_{CP^1}$
follow directly from  validity of the
Virasoro constraints.

{\noindent \it Remark 2}. In \cite{oko3} A. Okounkov and R. Pandharipande
proved that the Gromov - Witten potential of the {\it equivariant} GW invariants
of $CP^1$ and their descendents is the logarithm of the tau function
of the 2D Toda hierarchy of K. Ueno and K. Takasaki \cite{UT}. The tau function of
\cite{oko3} depends on an additional small parameter $t$. 
The non-equivariant limit corresponds to $t\to 0$.
It would be interesting
to derive the Lax representation of the extended Toda lattice by applying
a suitable limiting procedure to that of the 2D Toda hierarchy of \cite{UT}.
An interesting construction, due to Getzler \cite{getz-eq} of a nontrivial
reduction of the 2D Toda hierarchy depending on the parameter $t$ (it was called
equivariant Toda lattice) could give a clue to such a limiting procedure. We
plan to study the relationships between 2D and extended Toda hierarchies
in a subsequent publication.

\def\thetheorem{A.\arabic{theorem}}
\def\theprop{A.\arabic{prop}}
\def\thelemma{A.\arabic{lemma}}
\def\thecor{A.\arabic{cor}}
\def\theexam{A.\arabic{exam}}
\def\theremark{A.\arabic{remark}}
\def\theequation{A.\arabic{equation}}

\appendix
\makeatletter
\renewcommand{\@seccntformat}[1]{{Appendix:}\hspace{-2.3cm}}
\makeatother
\renewcommand{\thesection}{Appendix:}
\section{\quad\qquad \ \ Genus zero Gromov - Witten potential of $CP^1$.}
To compute the genus zero Gromov - Witten potential ${\cal F}_0({\bf t})$
according to the general scheme of \cite{npb, D1} one is to perform the
following computations (cf. \cite{DZ3}).

1. Compute the functions $\theta_{\al,p}(v,u)$ as the coefficients of expansion
of the following series
\eqa
&&\theta_1(v,u;z)=\sum_{p\geq 0} \theta_{1,p}(v,u) z^p
\nn\\
&&
=
-2\,e^{z v}\left(K_0(2 z e^{\frac12\,u})
+(\log{z}+\gamma) I_0(2 z e^{\frac12\,u})\right)\nn\\
&&=-2 e^{z\, v} \sum_{m\geq 0} (\gamma-{1\over 2} u + \psi(m+1)) e^{m\, u}
{z^{2\, m}\over (m!)^2},\label{theta-toda-1}\\
&&\theta_2(v,u;z)=\sum_{p\geq 0} \theta_{2,p}(v,u) z^p
\nn\\
&&
=
z^{-1}\,e^{z v}\,I_0(2 z e^{\frac12\,u})-z^{-1}
=z^{-1} \left(\sum_{m\geq 0}e^{m\, u + z \, v} {z^{2\, m}\over
(m!)^2}-1\right).\label{theta-toda-2}
\eeqa
Here $\gamma$ denotes Euler's constant, $\psi(z)$ stands for the digamma
function, $K_0(x)$ and $I_0(x)$ are modified Bessel functions.

2. Compute the functions $\Omega_{\al,p; \beta,q}^{[0]}(v,u)$ as the
coefficients of the following
generating series
\beq\label{matr}
\sum_{p,q\geq 0}\Omega_{\al,p; \beta,q}^{[0]}(v,u) z^p w^q
={1\over z+w} \left[ {\pal \theta_\al(v,u;z)\over\pal v}  
{\pal \theta_\beta(v,u;w)\over\pal u}+
{\pal \theta_\al(v,u;z)\over\pal u}  
{\pal \theta_\beta(v,u;w)\over\pal v}-\eta_{\al\beta}\right].
\eeq

3. Define the functions $v({\bf t})$, $u({\bf t})$ as the unique solution of the
system
\eqa\label{implicit}
&&
v=\sum t^{\beta,q} {\pal \theta_{\beta,q}\over \pal u}
\nn\\
&&
u=\sum t^{\beta,q} {\pal \theta_{\beta,q}\over \pal v}
\eeqa
having the expansion
\eqa
&&
v({\bf t}) = t^{1,0}+o(t)
\nn\\
&&
u({\bf t}) = t^{2,0}+o(t).
\nn
\eeqa

4. The genus zero Gromov - Witten potential of $CP^1$ is given by
\beq
{\cal F}_0({\bf t}) =\frac12 \sum \tilde t^{\al,p} \tilde t^{\beta,q}
\Omega_{\al,p; \beta,q}^{[0]}(v({\bf t}),u({\bf t})).
\eeq
Here
\beq
\tilde t^{\al,p}=\left\{ \matrix{ t^{1,1}-1, & \quad\al=1, \, p=1\cr
t^{\al,p}, & ~{\rm otherwise}.\cr}\right. 
\eeq
Let us write first few terms of the expansion of 
the resulting genus zero Gromov - Witten potential. For simplicity we denote
$$
t_p:= t^{1,p}, \quad s_p:= t^{2,p}, \quad p\geq 0.
$$
The potential is expanded in powers of $t_p$, $s_p$ and in $e^{s_0}$. The powers
of the exponential separate the intersection numbers on the moduli spaces 
$\bar M_{0,m}(CP^1,\beta)$ for different $\beta$. Thus, the terms without
exponential correspond to $\beta=0$ etc. In our expansion we collect the
intersection numbers up to degree $\beta=3$ and up to $m\leq 4$ punctures.
We also restrict the potential onto the subspace of couplings $t^{\al,p}$ with
$p\leq 3$. The maximal number of descendents is restricted to three.
\small{
\eqa
&&{\cal F}_0=\frac{t_0^2\,s_0}{2!}
+
\frac{t_0^2\,t_1\,s_0} {2!}+ \frac{t_0^3\,s_1}{3!}
+
\frac{t_0^3\,t_2\,s_0} {3!}
+ \frac{t_0^4\,s_2}{4!}
+\frac{t_0^2\,t_1^2\,s_0}{2!}
\nn\\
&&
+ 2\frac{t_0^3\,t_1\,s_1}{3!}
+
 \frac{t_0^4\,t_3\,s_0}{4!}
  +\frac {t_0^5\,s_3}{5!}
 + 3\frac{t_0^3\,t_1\,t_2\,s_0} {3!}
  + 3\frac{t_0^4\,t_2\,s_1}{4!}
  +3\frac{t_0^4\,t_1\,s_2} {4!}
\nn\\
&&  
+e^{s_0}\left[ 1 -2 \,t_1 +2\frac{ t_1^2}{2!} 
-2 \, t_0 t_2
+t_1s_0 +t_0 s_1
-
2\frac{{t_0}^2\,t_3}{2!}  
\right.
\nn\\
&&
 - 2\frac{{t_1}^2\,s_0}{2!} + 
  t_0\,t_2\,s_0 + \frac{{t_0}^2\,s_2}{2!}
-  
2\frac{{t_0}^2\,t_1\,t_3}{2!}   - 
  t_0\,t_1\,t_2\,s_0 + \frac{{t_0}^2\,t_3\,s_0}{2!}
  \nn\\
  &&
\left.   + 
  2\frac{{{t_1}^2\,{s_0}^2}}{(2!)^2} 
  - \frac{{t_0}^2\,t_2\,s_1}{2!}+ 
  t_0\,t_1\,s_0\,s_1 + \frac{{t_0}^2\,{s_1}^2}{2!} + 
  \frac{{t_0}^2\,t_1\,s_2}{2!} + \frac{{t_0}^3\,s_3}{3!}\right]
\nn\\
&&
+e^{2 s_0}\left[ -\frac34 t_3 +\frac14 s_2
+
\frac54\frac{{t_2}^2}{2!} + \frac34 t_1\,t_3 + 
  \frac14 t_3\,s_0 - \frac34 t_2\,s_1 + 
  \frac12 \frac{{s_1}^2}{2!} 
  \right.
  \nn\\
  &&
 - \frac14 {t_1\,s_2}
  + \frac14{t_0\,s_3}
+
2\,t_0\,t_2\,t_3 - \frac32\frac{{t_2}^2\,s_0}{2!} - 
  \frac32\,t_1\,t_3\,s_0 - 2\,t_0\,t_3\,s_1 
  \nn\\
  &&
  + 
  \frac12{t_2\,s_0\,s_1}- t_0\,t_2\,s_2 + 
  \frac12{t_1\,s_0\,s_2}+ t_0\,s_1\,s_2
+
2\,t_0\,t_1\,t_2\,t_3 + 4\frac{{t_0}^2\,{t_3}^2}{(2!)^2}
\nn\\
&&
  + 
  \frac{t_1\,{t_2}^2\,s_0}{2!} - 3\,t_0\,t_2\,t_3\,s_0 + 
  \frac{{t_2}^2\,{s_0}^2}{(2!)^2} + \frac{t_1\,t_3\,{s_0}^2}{2!} + 
  \frac{t_0\,{t_2}^2\,s_1}{2!} 
  \nn\\
  &&
  - 2\,t_0\,t_1\,t_3\,s_1 - 
  t_1\,t_2\,s_0\,s_1 + t_0\,t_3\,s_0\,s_1 - 
 2 \frac{ t_0\,t_2\,{s_1}^2 }{2!}
 + \frac{t_1\,s_0\,{s_1}^2}{2!} 
 \nn\\
 &&+ 
  3\frac{t_0\,{s_1}^3}{3!} - t_0\,t_1\,t_2\,s_2 - 
  3\frac{{t_0}^2\,t_3\,s_2}{2!} + t_0\,t_2\,s_0\,s_2 + 
  t_0\,t_1\,s_1\,s_2 
  \nn\\
  &&
\left.  + 2\frac{{t_0}^2\,{s_2}^2}{(2!)^2} - 
  \frac{{t_0}^2\,t_2\,s_3}{2!} + 
  \frac{t_0\,t_1\,s_0\,s_3}{2} + 
  \frac32\frac{{t_0}^2\,s_1\,s_3}{2!}\right]
 \nn\\
 &&
 +e^{3 s_0} \left[ \frac{50}{27}\frac{{t_3}^2}{2!} 
 - \frac79\,t_3\,s_2 + 
  \frac13\frac{{s_2}^2}{2!} - \frac29\,t_2\,s_3 + 
  \frac16{s_1\,s_3}
-
 2\frac{{t_2}^2\,t_3}{2!}  - \frac{14}9\frac{{t_3}^2\,s_0}{2!}
 \right.
 \nn\\
 &&
 + 
  2\,t_2\,t_3\,s_1 - 2\frac{t_3\,{s_1}^2}{2!} + 
  \frac{{t_2}^2\,s_2}{2!} + \frac13{t_3\,s_0\,s_2} - 
  t_2\,s_1\,s_2 + \frac{{s_1}^2\,s_2}{2!} - 
  t_0\,t_3\,s_3 
  \nn\\
  &&
 + \frac16{t_2\,s_0\,s_3}+ 
  \frac12{t_0\,s_2\,s_3}
-  
  4\frac{t_0\,t_2\,{t_3}^2}{2!} + 5\frac{{t_2}^2\,t_3\,s_0}{2!} + 
  4\frac{t_1\,{t_3}^2\,s_0}{2!} 
  + \frac23\frac{{t_3}^2\,{s_0}^2}{(2!)^2} 
  \nn\\
  &&
   + 
  8\frac{t_0\,{t_3}^2\,s_1}{2!} - 3\,t_2\,t_3\,s_0\,s_1 + 
  \frac{t_3\,s_0\,{s_1}^2}{2!} + 3\,t_0\,t_2\,t_3\,s_2 - 
  2\frac{{t_2}^2\,s_0\,s_2}{2!}
  \nn\\
  &&
  - 2\,t_1\,t_3\,s_0\,s_2 - 
  5\,t_0\,t_3\,s_1\,s_2 + t_2\,s_0\,s_1\,s_2 - 
 2\frac{ t_0\,t_2\,{s_2}^2}{2!} 
 + \frac{t_1\,s_0\,{s_2}^2}{2!} 
 \nn\\
 &&
 + 
 3 \frac{t_0\,s_1\,{s_2}^2}{2!} + 
  \frac{t_0\,{t_2}^2\,s_3}{2!} - t_0\,t_1\,t_3\,s_3 - 
  \frac12{t_1\,t_2\,s_0\,s_3} + 
  \frac12{t_0\,t_3\,s_0\,s_3} 
  \nn\\
  &&
 \left.   - 
  \frac32{t_0\,t_2\,s_1\,s_3}+ 
  \frac12{t_1\,s_0\,s_1\,s_3} +
  2\frac{t_0\,{s_1}^2\,s_3}{2!} + 
  \frac12{t_0\,t_1\,s_2\,s_3}
  + \frac{{t_0}^2\,{s_3}^2}{(2!)^2}\right].
  \nn
\eeqa
}


\end{document}